\documentclass[11pt]{amsart}

\usepackage{amsfonts,amssymb,amscd}
\newcommand{\R}{\mathbb R}
\newcommand{\Y}{\mathbb Y}

\newcommand{\A}{\mathbb A}
\newcommand{\B}{\mathbb B}
\newcommand{\Z}{\mathbb Z}
\newcommand{\N}{\mathbb N}
\newcommand{\te}{\theta}
\newcommand{\sm}{\sigma}
\newcommand{\si}{\sigma_k^{-1}}

\newcommand{\s}{\sim}
\newcommand{\bd}{\partial}
\newcommand{\st}{\stackrel}

\newcommand{\al}{\alpha}
\newcommand{\ph}{\varphi}

\newcommand{\Ga}{\Gamma}
\newcommand{\bu}{\bullet}
\newcommand{\ra}{\rightarrow}
\newcommand{\vl}{\; | \;}

\newtheorem{theo}{Theorem}

\newtheorem{lem}{Lemma}
\newtheorem{claim}{Claim}

\newcommand{\es}{\emptyset}
\title{ Three-page embeddings of singular knots}
\author[kurlin]{V.~kurlin$^*$}
\address{Department of Mathematics, Moscow State University,
Moscow, 119992, Russia}
\email {vak26@yahoo.com, kurlin@mccme.ru}
\thanks{$^*$
The first author was supported in part by grant INTAS YS 2001/2-30.}

\author[Vershinin]{V.~Vershinin$^{**}$}
\address{D\'epartement des Sciences Math\'ematiques,
Universit\'e Montpellier II,
Place Eug\'ene Bataillon,
34095 Montpellier cedex 5, France}
\email{ vershini@math.univ-montp2.fr}
\address{ Sobolev Institute of Mathematics, Novosibirsk, 630090,
Russia }
\email{ versh@math.nsc.ru}
\thanks{$^{**}$
The second author was supported in part by the French-Russian
Program of Research EGIDE (dossier No 04495UL)}

\begin{document}
\begin{abstract}
Construction of a semigroup with 15 generators and 84 relations
 is given. The center of this semigroup is in  one-to-one
 correspondence with the set of all isotopy classes of non-oriented
 singular knots (links with finitely many double intersections in general
 position) in $\R^3$.
\end{abstract}
\maketitle

\section{Introduction}

\subsection{Statement of the problem and results}
We develop the Dynnikov method of three-page embeddings for
 {\it links with singularities} of the following type: finitely many
double intersections in general position are possible.
More precisely, an algebraic solution of isotopy classification problem
 for non-oriented singular knots in $\R^3$ is given.
The key idea is a construction of a 3-page
 embedding for a neighborhood of a singular point.
In particular this construction gave a possibility of diminishing
the number of generators and defining relations in the semigroup
for singular knots.

\subsection{Review of the previous results}
An embedding of a link in a structure which looks like an open
 book with finitely many pages probably was considered for
 the first time by H.~Brunn in 1898 \cite{Br}.
Namely, he proved that each knot is isotopic to a knot that
 projects to a plane with only one singular point.
Later such embeddings were studied in papers \cite{BC1, BC2, Cr} and
 were used in \cite{BM}.
Moreover, these investigations provided a new link invariant --
\emph{arc index} \cite{BP, MB, Nu}.
It turned out that each link embeds in a book with only 3 pages.
In 1999 Dynnikov  classified all non-oriented links in $\R^3$
 up to an ambient isotopy encoding them by
 \emph{three-page diagrams} \cite{Dy1, Dy2}.
To be more precise we call these diagrams \emph{3-page embeddings}
 (see the definition in Subsection 2.1).
Dynnikov constructed some semigroup $DS$ such that there is a one-to-one
 correspondence between the center of $DS$ and the set of all isotopy classes of
 non-oriented links in $\R^3$.
Applying embeddings into a book with an arbitrary number of pages
 Dynnikov decreased the number of the relations in his semigroup $DS$
 \cite{Dy3}.
Analogously the first author obtained an isotopic
 classification of non-oriented knotted 3-valent graphs in $\R^3$ \cite{Ku}.

\subsection{Motivations}
In \cite{JM} singular knots were called
 \emph{chimerical graphs} and in \cite{Ka}
 four-valent graphs \emph{with rigid vertices}.
The recent study of singular knots and braids was motivated by the
 theory of Vassiliev invariants \cite{Bi93}.
The corresponding algebraic object is called the
{\it Baez-Birman monoid} or
 {\it singular braid monoid} \cite{Ba, Bi93}.
Some of its algebraic properties were investigated in \cite{DG, FKR, Ge2}.
For singular braids the analogue of Markov's theorem was
 proved \cite{Ge1}.
Many invariants of regular (non-singular) links, in particular, the
Alexander-Conway and Jones polynomials and Vassiliev invariants
 are extended to singular knots \cite{JM, Ka, KV, St}.
Homological properties of singular braids on
 infinitely many strings were studied by the second author \cite{Ve1, Ve2}.

\subsection{Basic definitions}
We work in the PL-category, i.e. images of circles
 under immersion in $\R^3$ are finite polygonal lines.
Formally, a \emph{ singular knot} is an immersion of several
 circles into $\R^3$ with possible double intersections
 in general position in a finite number of \emph{singular points}
(Fig.~1,~2).
Two \emph{branches} of a given singular knot pass through each
 singular point.
In the present paper (except Subsection~3.5) we consider only
 non-oriented singular knots, may be non-connected.
Note that a singular knot is a 4-valent graph embedded in $\R^3$.
An \emph{ambient PL-isotopy} between two graphs is a continuous
family of PL-homeomorphisms $\phi_t : \R^3 \to \R^3$, $t\in [0,1]$,
such that $\phi_0=\operatorname{id}$ and $\phi_1$ sends one of the
graphs to another. Singular knots are considered up to ambient
PL-isotopy that respects  the rigidity (or {\it template} by the
terminology of the work \cite{JM}) of each singular point.
If we omit the last restriction on isotopy then we come to the notion
of \emph{knuckle} 4-valent graph (graph with non-rigid vertices
 according to the terminology of the work \cite{Ka}).
Contrary to the case of singular knots an isotopy of knuckle graphs
can transpose edges at each 4-valent vertex.
In Subsection~3.5 we formulate classification Theorem~5 for knuckle graphs.
The same as for regular links, one can represent singular
 knots by plane diagrams equivalent up to the Reidemeister moves
 $R1-R5$ of Fig.~1 \cite{JM}.
We depict only PL-analogues of the corresponding smooth moves and
 omit subdivisions and extra breaks of edges.
In the case of knuckle 4-valent graphs the move $R5'$ is taken
instead of $R5$.
\bigskip

\begin{picture}(450,65)(0,0)

\put(0,45){\line(0,1){10}}
\put(0,5){\line(0,1){10}}
\put(15,15){\line(0,1){30}}
\put(0,45){\line(1,-2){15}}
\put(0,15){\line(1,2){5}}
\put(10,35){\line(1,2){5}}
\put(25,30){$\sim$}
\put(40,5){\line(0,1){50}}
\put(50,30){$\sim$}
\put(80,45){\line(0,1){10}}
\put(80,5){\line(0,1){10}}
\put(65,15){\line(0,1){30}}
\put(65,15){\line(1,2){15}}
\put(80,15){\line(-1,2){5}}
\put(70,35){\line(-1,2){5}}
\put(35,-10){$R1$}

\put(110,60){\line(1,-2){15}}
\put(125,30){\line(-1,-2){15}}
\put(125,60){\line(-1,-2){5}}
\put(110,30){\line(1,2){5}}
\put(110,30){\line(1,-2){5}}
\put(120,10){\line(1,-2){5}}
\put(135,30){$\sim$}
\put(150,5){\line(0,1){50}}
\put(160,5){\line(0,1){50}}
\put(170,30){$\sim$}
\put(185,30){\line(1,2){15}}
\put(185,30){\line(1,-2){15}}
\put(185,60){\line(1,-2){5}}
\put(195,40){\line(1,-2){5}}
\put(200,30){\line(-1,-2){5}}
\put(185,0){\line(1,2){5}}
\put(150,-10){$R2$}

\put(230,60){\line(1,-2){5}}
\put(240,40){\line(1,-2){15}}
\put(240,10){\line(1,2){5}}
\put(250,30){\line(1,2){5}}
\put(260,50){\line(1,2){5}}
\put(230,45){\line(1,0){35}}
\put(275,35){$\sim$}
\put(290,25){\line(1,0){35}}
\put(300,60){\line(1,-2){15}}
\put(320,20){\line(1,-2){5}}
\put(290,10){\line(1,2){5}}
\put(300,30){\line(1,2){5}}
\put(310,50){\line(1,2){5}}

\put(350,60){\line(1,-2){25}}
\put(360,10){\line(1,2){5}}
\put(370,30){\line(1,2){15}}
\put(362,45){\line(1,0){11}}
\put(345,45){\line(1,0){8}}
\put(382,45){\line(1,0){8}}
\put(395,35){$\sim$}
\put(422,25){\line(1,0){11}}
\put(405,25){\line(1,0){8}}
\put(442,25){\line(1,0){8}}
\put(420,60){\line(1,-2){25}}
\put(410,10){\line(1,2){15}}
\put(430,50){\line(1,2){5}}
\put(270,-10){$R3$}
\put(390,-10){$R3$}
\end{picture}
\medskip


\begin{picture}(450,80)(10,0)
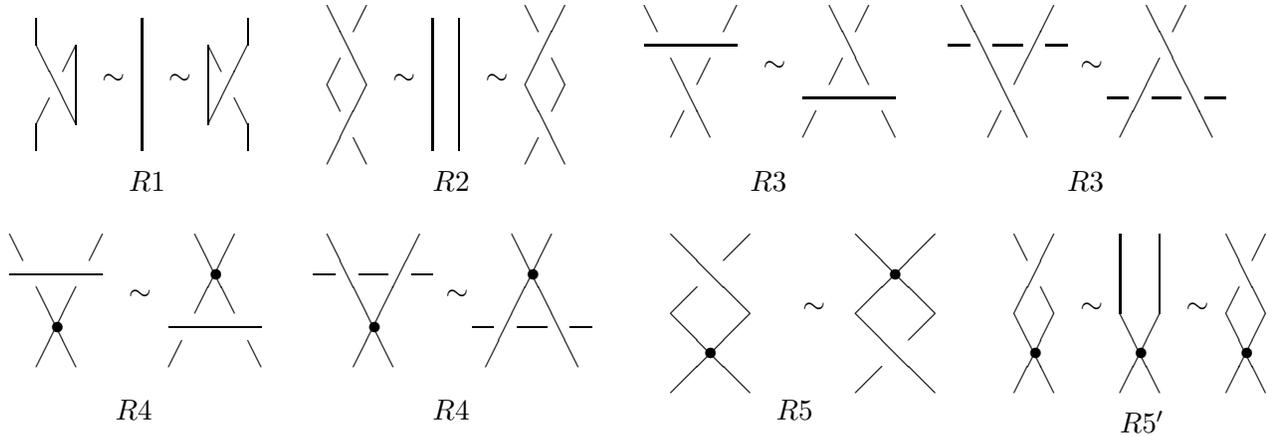
\put(80,-30){{\bf Fig.~1:} Reidemeister moves for singular knots
 and knuckle graphs}

\put(0,60){\line(1,-2){5}}
\put(10,40){\line(1,-2){15}}
\put(10,10){\line(1,2){15}}
\put(30,50){\line(1,2){5}}
\put(0,45){\line(1,0){35}}
\put(45,35){$\sim$}
\put(60,25){\line(1,0){35}}
\put(70,60){\line(1,-2){15}}
\put(90,20){\line(1,-2){5}}
\put(60,10){\line(1,2){5}}
\put(70,30){\line(1,2){15}}
\put(18,25){\circle*{4}}
\put(78,45){\circle*{4}}

\put(120,60){\line(1,-2){25}}
\put(130,10){\line(1,2){25}}
\put(132,45){\line(1,0){11}}
\put(115,45){\line(1,0){8}}
\put(152,45){\line(1,0){8}}
\put(165,35){$\sim$}
\put(192,25){\line(1,0){11}}
\put(175,25){\line(1,0){8}}
\put(212,25){\line(1,0){8}}
\put(190,60){\line(1,-2){25}}
\put(180,10){\line(1,2){25}}
\put(138,25){\circle*{4}}
\put(198,45){\circle*{4}}
\put(40,-10){$R4$}
\put(160,-10){$R4$}

\put(265,15){\circle*{4}}
\put(250,0){\line(1,1){30}}
\put(250,30){\line(1,1){10}}
\put(270,50){\line(1,1){10}}
\put(280,0){\line(-1,1){30}}
\put(280,30){\line(-1,1){30}}
\put(300,30){$\sim$}
\put(335,45){\circle*{4}}
\put(320,0){\line(1,1){10}}
\put(320,30){\line(1,1){30}}
 \put(340,20){\line(1,1){10}}
\put(350,0){\line(-1,1){30}}
\put(350,30){\line(-1,1){30}}
\put(290,-10){$R5$}

\put(388,15){\circle*{4}}
\put(380,0){\line(1,2){15}}
\put(395,0){\line(-1,2){15}}
\put(380,30){\line(1,2){15}}
\put(395,30){\line(-1,2){5}}
\put(385,50){\line(-1,2){5}}
\put(405,30){$\sim$}
\put(428,15){\circle*{4}}
\put(420,0){\line(1,2){15}}
\put(435,0){\line(-1,2){15}}
\put(420,30){\line(0,1){30}}
\put(435,30){\line(0,1){30}}
\put(445,30){$\sim$}
\put(468,15){\circle*{4}}
\put(460,0){\line(1,2){15}}
\put(475,0){\line(-1,2){15}}
\put(475,30){\line(-1,2){15}}
\put(460,30){\line(1,2){5}}
\put(470,50){\line(1,2){5}}
\put(420,-15){$R5'$}

\end{picture}
\vspace{1cm}

\subsection{The universal semigroup for singular knots.}
Everywhere the index $i$ belongs to the group $\Z_3=\{0,1,2\}$.
Consider the alphabet
 $\A=\lbrace a_i,b_i,c_i,d_i,x_i \vl i\in \Z_3 \rbrace$ with 15 letters
 (see Fig.~3 in Subsection~2.3 for their geometrical interpretation).
Let $SK$ be the semigroup on 15 generators of the alphabet $\A$
 and relations (1)-(10), which correspond to
 some "elementary ambient isotopies" of singular knots in $\R^3$.
$$\begin{array}{l}
(1) \quad a_i = a_{i+1} d_{i-1}, \quad b_i = a_{i-1} c_{i+1}, \quad
           c_i = b_{i-1} c_{i+1}, \quad d_i = a_{i+1} c_{i-1}, \\

(2) \quad x_i = d_{i+1}x_{i-1}b_{i+1}, \\

(3) \quad d_0 d_1 d_2 = 1,  \\

(4) \quad b_i d_i = d_i b_i =1, \\

(5) \quad d_i x_i d_i = a_i (d_i x_i d_i) c_i, \quad
          b_i x_i b_i = a_i (b_i x_i b_i) c_i, \\

(6) \quad x_i (d_{i+1} d_i d_{i-1}) = (d_{i+1} d_i d_{i-1}) x_i, \\

(7) \quad (d_{i} c_{i}) w = w (d_{i} c_{i}), \mbox{ where }
     w\in \lbrace c_{i+1}, x_{i+1}, b_{i} d_{i+1} d_{i} \rbrace, \\

(8) \quad (a_{i} b_{i}) w = w (a_{i} b_{i}), \mbox{ where }
     w\in \lbrace a_{i+1}, b_{i+1}, c_{i+1}, x_{i+1}, b_{i} d_{i+1}
d_{i} \rbrace, \\

(9) \quad t_i w = w t_i, \mbox{ where }
     t_i = b_{i+1} d_{i-1} d_{i+1} b_{i-1} ,
     w\in \lbrace a_i, b_i, c_i, x_i, b_{i-1} d_i d_{i-1} \rbrace, \\

(10) \quad (d_i x_i b_i) w = w (d_i x_i b_i), \mbox{ where }
     w\in \lbrace a_{i+1}, b_{i+1}, c_{i+1}, x_{i+1}, b_{i} d_{i+1}
d_{i} \rbrace.
  \end{array} $$
One relation in (4) is superfluous: it can be obtained from (3)
and the rest of the relations in (4).
Hence the total number of relations (1)-(10) is 84.

\subsection{Algebraic classification of singular knots}

\begin{theo}
Each singular knot can be represented by an element of the
semigroup $SK$.
\end{theo}
\begin{theo}
Two singular knots are ambiently isotopic in $\R^3$ if and only if
the corresponding elements of the semigroup $SK$ are equal.
\end{theo}
\begin{theo}
An arbitrary element of semigroup $SK$ corresponds to a singular knot
if and only if this element is central, i.e. it commutes with every
element of $SK$.
\end{theo}

As it will be shown in Theorem~4 the whole semigroup
 $SK$ describes a wider class of \emph{3-page singular tangles}.
The subsemigroup in $SK$ that is generated by 12 letters
 $a_i,b_i,c_i,d_i$ ($i\in \Z_3$) and 48 relations from (1)-(10)
 not containing letters $x_i$ ($i\in \Z_3$) coincides
 with Dynnikov's semigroup $DS$ of \cite{Dy1, Dy2}.
The center of semigroup $DS$ classifies all non-oriented
 (regular) links in $\R^3$ up to ambient isotopy.

\subsection{The content of the paper.}
In Subsection~2.1 we define three-page embeddings of singular knots.
Such embeddings are constructed from plane diagrams of knots
 in Subsection~2.2 and are encoded in Subsection~2.3.
Theorem~1 is proved in Subsection~2.5.
The ordinary singular tangles and three-page singular tangles
 are introduced in Subsection~3.1 and Subsection 3.2, respectively.
The later notion generalizes 3-page embeddings of singular knots
 and helps us in the proof of Theorem~2.
In Subsection~3.3 the three-page tangles are classified (Theorem~4).
Then Theorem~2 follows from Theorem~4 as a particular case.
Theorem~2 is applied in the proof of Theorem~3 in Subsection~3.4.
In Subsection~3.5 classification Theorem~5 for
 knuckle 4-valent graphs is formulated.
In Section~4 we deduce Lemma 3 which is used in the proof of Theorem~4.

\subsection{Acknowledgments}
The first author is thankful for hospitality to  I.~K.~Babenko,
 J.~Lafontaine and to the University Montpellier~II (France) where
 this paper was written.
Also he thanks I.~A.~Dynnikov for his interest in our results and
 the scientific advisor V.~M.~Buchstaber for support.
The authors are thankful to S.~S.~Kutateladze for help and advises in
the presentation of the paper.

\section {Three-page embeddings}

\subsection{Formal definition of three-page embeddings}
An \emph{ arc} of a singular knot $K\subset \R^3$ with endpoint
$A\in K$ is a sufficiently small segment $J\subset K$ with
$A\in\partial J$.
Then 4 arcs issue from each singular point.
Let $P_0, P_1$ and $P_2$ be three half-planes in $\R^3$ with a common
 oriented boundary:
 $\partial P_0=\partial P_1=\partial P_2=\alpha$ (Fig.~2).
Put $\Y = P_0\cup P_1\cup P_2$ and call this union a
 \emph{book with three pages}.
An embedding of a singular knot $K$ in book $\Y$ is called a
 \emph{3-page embedding}, if the following conditions hold:

 1) all singular points of $K$ lie on the axis $\alpha$;

 2) \emph{finiteness}: the intersection
  $K \cap \alpha = A_1\cup \dots \cup A_m$ is a finite point set;

 3) at every non-singular point $A_j\in K\cap\al$ two arcs lie in
different half-planes;

 4) a neighborhood of a singular point $A_j$ lies in
  the plane $P_{i-1}\cup P_{i+1}$ for some $i\in\Z_3$;

 5) \emph{monotonicity}: for each $i\in \Z_3$ the restriction of
  the orthogonal projection $\R^3\to \alpha \approx \R$ to
  each connected component of $K\cap P_i$ is a monotone function.

\subsection{Construction of a 3-page embedding from a plane diagram}
Let $D$ be a plane diagram of a singular knot $K$, i.e. a planar
4-valent graph with vertices of two types: one corresponds to singular
points of $K$ and the other denotes the usual \emph{crossings} in
a planar representation of $K$.
Given singular point $B$, let us mark two arcs with endpoint $B$,
namely a \emph{singular bridge} $L_B$ lying on
 different branches of our singular knot.
Also given crossing of the diagram $D$, mark a small segment
 (a \emph{regular bridge}) in the overcrossing arc (Fig.~2).

\begin{picture}(480,120)(0,-5)

\put(10,90){\vector(1,-1){10}}
\put(20,80){\line(1,-1){10}}
\put(50,70){\vector(1,-2){10}}
\put(60,50){\line(1,-2){10}}
\put(80,40){\vector(1,1){10}}
\put(80,60){\line(1,-1){10}}
\put(80,80){\vector(1,1){20}}
\put(100,80){$\alpha$}
\put(40,-10){$D$}
\put(85,20){$B$}
\put(67,39){{\footnotesize $L_B$}}
\put(105,35){$L_2$}
\put(75,0){$L_1$}
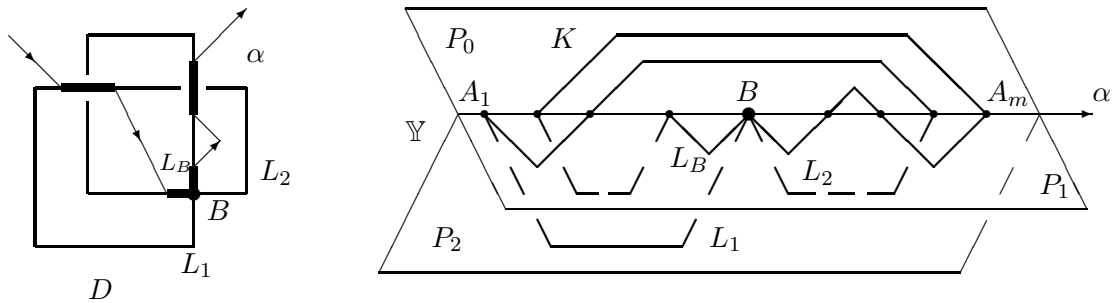
\put(80,-20){{\bf Fig.~2:} Three-page embedding $K\subset\Y$,
 $w_K=a_0 a_1 b_2 b_0 x_0 b_2 d_2 c_1 c_2$.}
\put(160,50){$\Y$}

\put(180,60){\vector(1,0){240}}
\put(420,65){$\alpha$}
\put(150,0){\line(1,0){220}}
\put(160,100){\line(1,0){220}}
\put(198,24){\line(1,0){220}}
\put(160,100){\line(1,-2){38}}
\put(380,100){\line(1,-2){38}}
\put(150,0){\line(1,2){30}}
\put(400,60){\line(-1,-2){5}}
\put(390,40){\line(-1,-2){5}}
\put(380,20){\line(-1,-2){10}}
\put(170,10){$P_2$}
\put(400,30){$P_1$}

\thicklines

\put(50,70){\line(1,0){25}}
\put(40,30){\line(1,0){60}}
\put(20,70){\line(1,0){55}}
\put(85,70){\line(1,0){15}}
\put(20,10){\line(1,0){60}}
\put(40,90){\line(1,0){40}}

\put(20,10){\line(0,1){60}}
\put(40,30){\line(0,1){35}}
\put(40,75){\line(0,1){15}}
\put(80,10){\line(0,1){80}}
\put(100,30){\line(0,1){40}}

\put(190,60){\line(1,-1){20}}
\put(190,60){\line(1,-2){5}}
\put(200,40){\line(1,-2){5}}
\put(210,20){\line(1,-2){5}}
\put(210,40){\line(1,1){40}}
\put(210,60){\line(1,1){30}}
\put(215,10){\line(1,0){50}}
\put(210,60){\line(1,-2){5}}
\put(220,40){\line(1,-2){5}}
\put(225,30){\line(1,0){8}}
\put(237,30){\line(1,0){8}}
\put(245,30){\line(1,2){5}}
\put(255,50){\line(1,2){5}}

\put(240,90){\line(1,0){110}}
\put(350,90){\line(1,-1){20}}

\put(250,80){\line(1,0){90}}

\put(260,60){\line(1,-1){15}}

\put(265,10){\line(1,2){5}}
\put(275,30){\line(1,2){5}}

\put(290,60){\line(1,-1){15}}
\put(290,60){\line(-1,-1){15}}
\put(290,60){\line(1,-2){5}}
\put(290,60){\line(-1,-2){5}}
\put(300,40){\line(1,-2){5}}
\put(305,30){\line(1,0){10}}
\put(320,30){\line(1,0){10}}
\put(335,30){\line(1,0){10}}

\put(305,45){\line(1,1){25}}
\put(330,70){\line(1,-1){30}}
\put(360,60){\line(-1,1){20}}
\put(360,60){\line(-1,-2){5}}
\put(350,40){\line(-1,-2){5}}

\put(380,60){\line(-1,-1){20}}
\put(380,60){\line(-1,1){30}}

\put(180,65){$A_1$}
\put(380,65){$A_m$}
\put(175,85){$P_0$}
\put(215,85){$K$}
\put(285,65){$B$}
\put(260,40){$L_B$}
\put(275,10){$L_1$}
\put(310,35){$L_2$}

\linethickness{1mm}
\put(30,70){\line(1,0){20}}
\put(80,60){\line(0,1){20}}
\put(70,30){\line(1,0){10}}
\put(80,30){\line(0,1){10}}

\put(80,30){\circle*{5}}
\put(190,60){\circle*{3}}
\put(210,60){\circle*{3}}
\put(230,60){\circle*{3}}
\put(260,60){\circle*{3}}
\put(290,60){\circle*{5}}
\put(320,60){\circle*{3}}
\put(340,60){\circle*{3}}
\put(360,60){\circle*{3}}
\put(380,60){\circle*{3}}

\end{picture}
\vspace{0.5cm}

\noindent
Then take a non-self-intersected oriented path $\alpha$ in the plane
 of the diagram $D$ with the following properties:

1) the endpoints of the path $\al$ lie far from $D$;

2) the path $\alpha$ traverses each bridge only once;

3) \emph{transversality}: the path $\alpha$ intersects the
  diagram of $D$ transversally beyond the bridges;

4) \emph{balance}: for a singular point $B$ consider
 two non-marked arcs $L_1, L_2$ with endpoint $B$, not containing
 the singular bridge $L_B$;
 then one of these arcs has to meet by the second endpoint the path
 $\alpha$ to the left of $L_B$ and the other has to meet by
 the second endpoint the path $\alpha$ to the right of $L_B$
(Fig.~2).
\smallskip

Such a path $\alpha$ can be easily found as follows:
consider only bridges in the plane, i.e. finitely many arcs.
Pass an arbitrary path $\alpha$ through each bridge satisfying
1) and 2). Then the transversality property 3) will hold, if we
move our path $\alpha$ in general position with respect to the
diagram $D$. Suppose that for the resulting path the balance
property 4) does not hold for a singular point $B$, i.e. both
non-marked arcs $L_1,L_2$ with endpoint $B$ meet by the second
endpoint the path $\al$ to the left of the bridge $L_B$
(for example). Then slightly move the path $\al$ to the right
of $L_B$ using a move like the Reidemeister move $R2$ such that
one of the two non-marked arcs $L_1,L_2$ (this is the arc $L_2$
in Fig.~2) meets by the second endpoint the path $\al$ to the
right of $L_B$.

Now deform the plane of $D$ in such a way that $\alpha$ becomes
a straight line and the following \emph{monotonicity} condition
holds~: the restriction of the orthogonal projection
$\R^2\to \alpha \approx \R$ to each connected component
of $D-\alpha$ is a monotonic function.
Denote by $P_0$ the upper half-plane over $\alpha$ and the
lower half-plane by $P_2$ (Fig.~2).
Finally, attach the third half-plane $P_1$ at $\alpha$ (at the
reader's side)
and push  out all bridges into $P_1$ according to the following rules.
Each regular bridge becomes a trivial arc.
Each singular bridge becomes a broken line, which looks like
 the letter ``W" that meets the axis $\alpha $ in its 3 upper vertices,
 and the middle one is the singular point $B$ (Fig.~2).
In fact, a neighborhood of any singular point can be embedded into
the plane $P_0\cup P_2$
 not pushing out marked arcs into the third half-plane $P_1$.
We used the notion of singular bridge for simplifying our argument.

\subsection{Encoding of three-page embeddings}
Each 3-page embedding is uniquely determined by its small neighborhood
 near the axis $\alpha$ in the book $\Y$.
Indeed, in order to reconstruct the whole embedding it is sufficient
 to connect the opposite-directed arcs in each half-plane starting from
 interior arcs.
We always mean that the half-plane $P_0$ is above the
 axis $\al$, and the half-planes $P_1,P_2$ are below $\al$.
Moreover, we suppose that $P_1$ is above $P_2$, i.e.
 arcs in $P_2$ are drawn by dashed lines.
Only the following 15 pictures may occur in a 3-page embedding of
a singular knot near the axis $\alpha$:

\begin{picture}(480,65)(0,15)

\put(0,60){\vector(1,0){60}}
\put(30,60){\circle*{3}}
\put(90,60){\vector(1,0){60}}
\put(120,60){\circle*{3}}
\put(180,60){\vector(1,0){60}}
\put(210,60){\circle*{3}}
\put(270,60){\vector(1,0){60}}
\put(300,60){\circle*{3}}
\put(360,60){\vector(1,0){60}}
\put(390,60){\circle*{5}}

\thicklines

\put(30,60){\line(1,-2){15}}
\put(30,60){\line(1,-1){10}}
\put(45,45){\line(1,-1){10}}
\put(20,20){$a_0$}

\put(120,60){\line(1,-2){15}}
\put(120,60){\line(-1,-1){10}}
\put(105,45){\line(-1,-1){10}}
\put(120,20){$b_0$}

\put(210,60){\line(-1,-2){15}}
\put(210,60){\line(-1,-1){10}}
\put(195,45){\line(-1,-1){10}}
\put(210,20){$c_0$}

\put(300,60){\line(-1,-2){15}}
\put(300,60){\line(1,-1){10}}
\put(315,45){\line(1,-1){10}}
\put(300,20){$d_0$}

\put(390,60){\line(1,-2){15}}
\put(390,60){\line(-1,-2){15}}
\put(390,60){\line(1,-1){10}}
\put(390,60){\line(-1,-1){10}}
\put(405,45){\line(1,-1){10}}
\put(375,45){\line(-1,-1){10}}
\put(390,20){$x_0$}

\end{picture}

\begin{picture}(480,80)(0,15)

\put(0,60){\vector(1,0){60}}
\put(30,60){\circle*{3}}
\put(90,60){\vector(1,0){60}}
\put(120,60){\circle*{3}}
\put(180,60){\vector(1,0){60}}
\put(210,60){\circle*{3}}
\put(270,60){\vector(1,0){60}}
\put(300,60){\circle*{3}}
\put(360,60){\vector(1,0){60}}
\put(390,60){\circle*{5}}
\thicklines

\put(30,60){\line(1,1){25}}
\put(30,60){\line(1,-1){10}}
\put(45,45){\line(1,-1){10}}
\put(30,20){$a_1$}

\put(120,60){\line(-1,1){25}}
\put(120,60){\line(1,-1){10}}
\put(135,45){\line(1,-1){10}}
\put(120,20){$b_1$}

\put(210,60){\line(-1,1){25}}
\put(210,60){\line(-1,-1){10}}
\put(195,45){\line(-1,-1){10}}
\put(210,20){$c_1$}

\put(300,60){\line(1,1){25}}
\put(300,60){\line(-1,-1){10}}
\put(285,45){\line(-1,-1){10}}
\put(300,20){$d_1$}

\put(390,60){\line(1,1){25}}
\put(390,60){\line(-1,1){25}}
\put(390,60){\line(-1,-1){10}}
\put(375,45){\line(-1,-1){10}}
\put(390,60){\line(1,-1){10}}
\put(405,45){\line(1,-1){10}}
\put(390,20){$x_1$}

\end{picture}

\begin{picture}(450,80)(0,15)
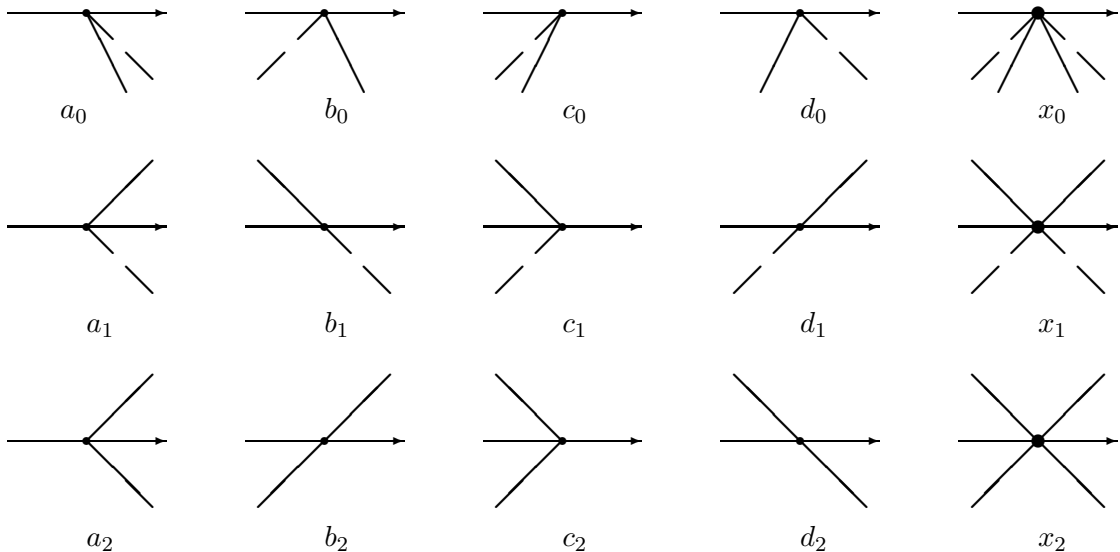
\put(80,-5){{\bf Fig.~3:} Geometric interpretation of letters of
the alphabet $\A$}

\put(0,60){\vector(1,0){60}}
\put(30,60){\circle*{3}}
\put(90,60){\vector(1,0){60}}
\put(120,60){\circle*{3}}
\put(180,60){\vector(1,0){60}}
\put(210,60){\circle*{3}}
\put(270,60){\vector(1,0){60}}
\put(300,60){\circle*{3}}
\put(360,60){\vector(1,0){60}}
\put(390,60){\circle*{5}}
\thicklines

\put(30,60){\line(1,1){25}}
\put(30,60){\line(1,-1){25}}
\put(30,20){$a_2$}

\put(120,60){\line(1,1){25}}
\put(120,60){\line(-1,-1){25}}
\put(120,20){$b_2$}

\put(210,60){\line(-1,1){25}}
\put(210,60){\line(-1,-1){25}}
\put(210,20){$c_2$}

\put(300,60){\line(-1,1){25}}
\put(300,60){\line(1,-1){25}}
\put(300,20){$d_2$}

\put(390,60){\line(-1,1){25}}
\put(390,60){\line(1,1){25}}
\put(390,60){\line(1,-1){25}}
\put(390,60){\line(-1,-1){25}}
\put(390,20){$x_2$}

\end{picture}
\vspace{1cm}

Let $W$ be the set of all words on the alphabet
$\A=\lbrace a_i, b_i, c_i, d_i, x_i \vl i\in \Z_3\rbrace$
 including the empty word $\varnothing$.
For a given 3-page embedding of the knot $K$ write one by one letters
of $\A$ corresponding to the intersection points of $K\cap \alpha$.
We obtain some word $w_K \in W$ (Fig.~2).

\subsection{Balanced words}
Note that by encoding of Subsection~2.3 one cannot obtain all the
words of $W$.
We call a word \emph{balanced}, if it encodes some 3-page embedding.
The following simple geometric criterion for a word to be balanced
is available~:
 in each half-plane $P_i$ all arcs are connected with each other.
Arcs of non-balanced 3-page embedding can go to infinity not
 meeting each other.
One can easily rewrite this criterion algebraically in terms of
 the alphabet $\A$.
For $i\in \Z_3$ a word $w$ is called \emph{$i$-balanced}, if
 after the following substitution
 $$a_i,b_i,c_i,d_i,x_i \ra \varnothing,  \quad
   a_{i\pm 1},b_{i-1},d_{i+1} \ra (,  \quad
   b_{i+1},c_{i\pm 1},d_{i-1} \ra ),  \quad
   x_{i\pm 1} \ra )(  $$
 we obtain an expression with \emph{completely balanced brackets}
 (or with \emph{correctly placed brackets} in another terminology).
This means that in each place the number of left brackets is not less
 than the number of the right ones, and their total numbers are equal.
By $W_i$ we denote the set of all $i$-balanced words in $\A$.
Then a word $w$ is called \emph{balanced}, if it is $i$-balanced for
 each $i\in \Z_3$.
So, the set of all balanced words is
 $W_b=W_0\cap W_1\cap W_2\subset W$.

\subsection{Proof of Theorem 1}
Take a plane diagram $D$ of a given singular knot $K$.
Starting with the diagram $D$ construct a 3-page embedding
$K\subset \Y$ described in Subsection 2.2.
Encode the obtained 3-page embedding of $K$ by the balanced
 word $w_K\in W_b$ according to the rules of Subsection~2.3.
Finally, consider the word $w_K$ as an element of the semigroup $SK$.
\qed

\section{Singular tangles}

\subsection{Semigroup $ST$ of singular tangles}
In order to prove Theorem~2 we need the notion of singular tangle.
The category of tangles (without singularities) was studied by
V.~G.~Turaev \cite{Tu}.
Take two horizontal semilines $\R_+\subset \R^3$, for
example given by coordinates: $(r, 0, 0)$ and $(r, 0, 1)$, where
$r\in \R_+$.
Mark the natural points $(j,0,0)$, $(k,0,1)$ for all $j,k\in \N$ on
both semilines.
Let $\Ga$ be a non-connected non-oriented infinite graph $\Ga$
 with vertices of valency 1 and 4.
A \emph{singular tangle} is an embedding of $\Ga$ into
 the 3-dimensional layer $\{0\leq z\leq 1\}$ such that (Fig.~4):

1) the set of the 1-valent vertices of the graph $\Ga$ coincides with
 the set of marked  points
 $$\{\; (j,0,0), \; (k,0,1) \vl j,k\in \N \;\};$$

2) all connected components of the graph $\Ga$ lying sufficiently far
 from the origin are the line segments connecting between the
 points $(k,0,0)$ and  $(j,0,1)$ such that the difference $k-j$
 is constant for all large $j$;

3) there exists a neighborhood of each 4-valent vertex of the graph
$\Ga$ which lies in a plane.
\smallskip

We consider singular tangles up to ambient isotopy in the layer
 $\{0\leq z\leq 1\}$ fixed on its boundary and such that condition~3)
holds. Singular tangles can be represented by their plane diagrams
analogous to singular knots (Fig.~4).
One can obtain a {\it product} of singular tangles
 $\Ga_1 \times \Ga_2$ by attaching the top semiline of $\Ga_2$ to
 the bottom semiline of $\Ga_1$.
So, the isotopy classes of singular tangles form some semigroup $ST$.
The \emph{unit} of $ST$ is the singular tangle consisting of
vertical segments.
Let us introduce the singular
 tangles: $\xi_k$, $\eta_k$, $\sm_k$, $\si$, $\tau_k$ ($k\in \N$):

\begin{picture}(360,75)(-30,-20)
\thicklines

\put(0,35){\line(1,-1){15}}
\put(15,20){\line(1,1){15}}
\put(-10,5){\line(0,1){30}}
\put(0,5){\line(3,1){30}}
\put(30,15){\line(1,2){10}}
\put(30,5){\line(1,1){10}}
\put(-5,40){\footnotesize $k$}
\put(15,40){\footnotesize $k+1$}
\put(10,-10){$\xi_k$}

\put(0,35){\circle*{3}}
\put(30,35){\circle*{3}}
\put(0,5){\circle*{3}}
\put(30,5){\circle*{3}}
\put(-10,5){\circle*{3}}
\put(-10,35){\circle*{3}}

\put(40,5){\circle*{3}}
\put(40,35){\circle*{3}}

\put(70,5){\line(0,1){30}}
\put(80,5){\line(1,1){15}}
\put(95,20){\line(1,-1){15}}
\put(80,35){\line(3,-1){30}}
\put(110,25){\line(1,-2){10}}
\put(110,35){\line(1,-1){10}}
\put(75,40){\footnotesize $k$}
\put(95,40){\footnotesize $k+1$}
\put(90,-10){$\eta_k$}

\put(70,5){\circle*{3}}
\put(70,35){\circle*{3}}
\put(80,35){\circle*{3}}
\put(110,35){\circle*{3}}
\put(80,5){\circle*{3}}
\put(110,5){\circle*{3}}
\put(120,5){\circle*{3}}
\put(120,35){\circle*{3}}

\put(150,5){\line(0,1){30}}
\put(160,5){\line(1,1){30}}
\put(160,35){\line(1,-1){12}}
\put(178,17){\line(1,-1){12}}
\put(200,5){\line(0,1){30}}
\put(155,40){\footnotesize $k$}
\put(175,40){\footnotesize $k+1$}

\put(150,5){\circle*{3}}
\put(150,35){\circle*{3}}
\put(160,35){\circle*{3}}
\put(190,35){\circle*{3}}
\put(160,5){\circle*{3}}
\put(190,5){\circle*{3}}
\put(200,5){\circle*{3}}
\put(200,35){\circle*{3}}
\put(170,-10){$\sm_k$}

\put(230,5){\line(0,1){30}}
\put(240,5){\line(1,1){12}}
\put(258,23){\line(1,1){12}}
\put(240,35){\line(1,-1){30}}
\put(280,5){\line(0,1){30}}
\put(235,40){\footnotesize $k$}
\put(255,40){\footnotesize $k+1$}

\put(230,5){\circle*{3}}
\put(230,35){\circle*{3}}
\put(240,35){\circle*{3}}
\put(240,5){\circle*{3}}
\put(270,5){\circle*{3}}
\put(270,35){\circle*{3}}
\put(280,5){\circle*{3}}
\put(280,35){\circle*{3}}
\put(250,-10){$\si$}

\put(310,5){\line(0,1){30}}
\put(320,5){\line(1,1){30}}
\put(320,35){\line(1,-1){30}}
\put(360,5){\line(0,1){30}}
\put(315,40){\footnotesize $k$}
\put(335,40){\footnotesize $k+1$}

\put(310,5){\circle*{3}}
\put(310,35){\circle*{3}}
\put(320,5){\circle*{3}}
\put(320,35){\circle*{3}}
\put(335,20){\circle*{5}}
\put(350,5){\circle*{3}}
\put(350,35){\circle*{3}}
\put(360,5){\circle*{3}}
\put(360,35){\circle*{3}}
\put(330,-10){$\tau_k$}
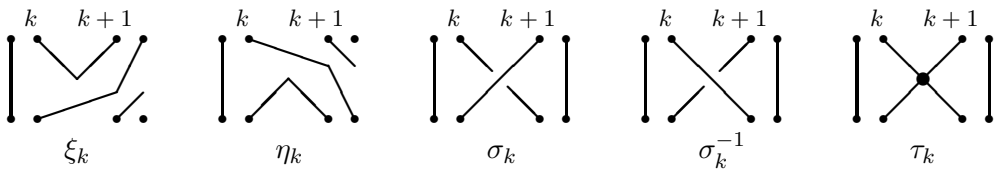
\put(50,-30){{\bf Fig.~4:} Generators of the singular tangles}

\end{picture}
\vspace{1cm}

The following lemma transfers results of \cite{Tu} from
 the classical case to ours.

\begin{lem}
The semigroup $ST$ of singular tangles is generated by the elements
 $\xi_k$, $\eta_k$, $\sm_k$, $\si$, $\tau_k$, $k\in \N$ (Fig.~4)
 and relations (11)-(23), where $k,l\in \N$ :
$$\begin{array}{llllll}
(11) & \xi_k \xi_l = \xi_{l+2} \xi_k, &
        \xi_k \eta_l = \eta_{l+2} \xi_k, &
        \xi_k \sm_l = \sm_{l+2} \xi_k, &
        \xi_k \tau_l = \tau_{l+2} \xi_k &
        (l\geq k); \\

(12) & \eta_k \xi_l = \xi_{l-2} \eta_k, &
        \eta_k \eta_l = \eta_{l-2} \eta_k, &
        \eta_k \sm_l = \sm_{l-2} \eta_k, &
        \eta_k \tau_l = \tau_{l-2} \eta_k &
        (l\geq k+2); \\

(13) & \sm_k \xi_l = \xi_{l} \sm_k, &
        \sm_k \eta_l = \eta_{l} \sm_k, &
        \sm_k \sm_l = \sm_{l} \sm_k, &
        \sm_k \tau_l = \tau_{l} \sm_k &
        (l\geq k+2); \\

(14) & \tau_k \xi_l = \xi_{l} \sm_k, &
        \tau_k \eta_l = \eta_{l} \sm_k, &
        \tau_k \sm_l = \sm_{l} \sm_k, &
        \tau_k \tau_l = \tau_{l} \tau_k &
        (l\geq k+2);
\end{array}$$

$$\begin{array}{ll}
(15) \quad \eta_{k+1} \xi_k = 1 = \eta_{k} \xi_{k+1}; &

(16) \quad \eta_{k+2} \sm_{k+1} \xi_k = \si =
           \eta_k \sm_{k+1} \xi_{k+2}; \\

(17)  \quad \eta_{k+2} \tau_{k+1} \xi_k = \tau_k =
           \eta_k \tau_{k+1} \xi_{k+2}; &

(18) \quad \eta_k \sm_k = \eta_k, \; \sm_k \xi_k = \xi_k; \\

(19) \quad \sm_k \si = 1 = \si \sm_k; &

(20) \quad \sm_k \sm_{k+1} \sm_k = \sm_{k+1} \sm_k \sm_{k+1}; \\

(21) \quad \sm_k \sm_{k+1} \tau_k = \tau_{k+1} \sm_k \sm_{k+1}; &

(22) \quad \tau_k \sm_{k+1} \sm_k = \sm_{k+1} \sm_k \tau_{k+1}; \\

(23) \quad  \sm_k \tau_k = \tau_k \sm_k. & \\
\end{array} $$
\end{lem}
\begin{proof}
Recall that we work in the $PL$-category.
It means that a given singular tangle
in the layer $\{0\leq z\leq 1\}$ consists of finite broken lines.
The local maxima and minima of the height function are called
\emph{extremal points}.
We call a \emph{peculiarity} of a diagram of a tangle either a
4-valent vertex, or a crossing, or an extremal point.
We say that a singular tangle is in \emph{general position} if
 its plane diagram satisfies the following conditions:

 1) the set of all peculiarities is finite;

 2) crossings do not coincide with extremal points;

 3) for each 4-valent vertex two arcs go up, and the rest two go down;

 4) each horizontal line (that is parallel to the $Ox$-axis)
contains at most than one peculiarity.

Obviously, by a slight deformation every tangle can be moved in a
general position. Then the tangle diagram is splitted by horizontal
lines into bands such that each of them contains only one peculiarity.
Considering the peculiarities from the top to the bottom one by one,
 write the corresponding generators from Fig.~4 from left to right.
Namely, the generators $\xi_k, \eta_k$ represent extremal points;
the generators $\sm_k, \si$ correspond to crossings; $\tau_k$ presents
a 4-valent vertex.
It suffices to show that every ambient isotopy of singular tangles
decomposes into ''elementary isotopies" corresponding to relations
(11)-(23).
It follows from the Reidemeister theorem \cite{JM} and the general
position arguments that an arbitrary isotopy of singular tangles can be
decomposed into the following moves:

 1) general position isotopy;

 2) swopping of heights of two peculiarities;

 3) creation or annihilation of a couple of close extremal points;

 4) an isotopy of a crossing or a 4-valent vertex near extremal point;

 5) the Reidemeister moves $R1-R5$ (Fig.~1).

The first type isotopies keep the constructed word in the letters
$\xi_k$ ,$\eta_k$, $\sm_k$, $\si$, $\tau_k$, $k\in \N$.
The second type isotopies are desribed by relations (11)-(14);
the third type isotopies correspond to relations (15).
In \cite[proof of lemma 3.4]{Tu} it was shown that in the smooth
category all isotopies of a crossing near extremal point
are geometrically decomposed into relations (16).
Similarly, in our $PL$-case we can check that relations
(17) are sufficient to isotope a 4-valent vertex near extremal point.
Finally, Reidemeister moves $R1-R5$ correspond to relations
(18)-(23), respectively.
\end{proof}

\subsection{Three-page singular tangles}
A notion of 3-page singular tangle will be used in the
 proofs of Theorems~2 and 3. Consider three semi-lines in the
horizontal plane $\{z=0\}$ having a common endpoint. Let it be
for example:
$$T=\{x\geq 0,y=z=0\}\cup \{y\geq 0,x=z=0\}\cup \{x\leq 0,y=z=0\}
\subset \{z=0\}.$$
Mark the integer points on the semilines:
 $\{(j,0,0),(0,k,0),(-l,0,0) \vl j,k,l\in \N \}$.
Let $I$ be an interval connecting the points $(0,0,0)$ and $(0,0,1)$.
Put:
 $$P_0=\{x\geq 0,y=z=0\}\times I, \quad
   P_1=\{y\geq 0,x=z=0\}\times I, \quad
   P_2=\{x\leq 0,y=z=0\}\times I.$$
Formally, here $P_i$ is not a half-plane, but a band $I\times\R$,
 which we call a \emph{page}.
The book $\Y$ of Subsection~2.1 is the interior of the set
$T\times I$, i.e. in Section~2 we considered the embeddings
$K\subset T\times I$ such that $K\cap \bd(T\times I)=\emptyset$.
Let $\Gamma$ be a non-connected non-oriented infinite graph with
vertices of valency 1 and 4.
A \emph{three-page singular tangle} is an embedding of $\Gamma$
 into a book $T\times I$ such that (Fig.~5):

 1) the set of 1-valent vertices of the graph $\Gamma$ coincides with
  the set of the marked points
 $$\{(j,0,0),(j,0,1), (0,k,0),(0,k,1), (-l,0,0),(-l,0,1),
\vl j,k,l\in \N \};$$

 2) all 4-valent vertices of $\Ga$ lie in the segment $I$;

 3) \emph{finiteness}: the intersection
  $\Gamma \cap I = A_1\cup \dots \cup A_m$ is a finite point set;

 4) the two arcs of any 2-valent vertex $A_j\in \Ga\cap I$ lie
in different half-planes;

 5) a neighborhood of each 4-valent vertex of $\Gamma$ lies
 in one pair of pages from $P_0,P_1,P_2$;

 6) \emph{monotonicity}: for every $i\in \Z_3$ restriction of the
orthogonal projection $T\times I\to I \approx [0,1]$ to each
connected component of $\Gamma \cap P_i$ is a monotone function.

 7) for each $i\in \Z_3$ all connected components of the graph
$\Gamma$ lying in a plane $P_i$ sufficiently far from the origin
are parallel line segments.
\vspace{0.5cm}

\begin{picture}(360,70)(-30,0)

\put(-10,70){\line(1,0){65}}
\put(20,70){\line(-2,-1){35}}
\put(15,30){\line(1,0){40}}
\put(-5,30){\line(1,0){10}}
\put(20,30){\line(-2,-1){35}}
\put(20,70){\vector(0,-1){55}}

{\thicklines
\put(30,70){\line(-1,-1){20}}
\put(20,40){\line(-1,1){10}}
\put(20,40){\line(2,3){20}}
\put(30,30){\line(1,2){20}}
\put(40,30){\line(1,2){10}}
\put(10,65){\line(-1,-3){5}}
\put(5,50){\line(1,-5){5}}
\put(0,30){\line(0,1){5}}
\put(0,40){\line(0,1){5}}
\put(0,50){\line(0,1){5}}
\put(0,65){\line(0,1){5}}
\put(-10,15){\line(0,1){40}}
}

\put(0,70){\circle*{3}}
\put(0,30){\circle*{3}}
\put(10,66){\circle*{3}}
\put(10,25){\circle*{3}}
\put(-10,55){\circle*{3}}
\put(-10,15){\circle*{3}}
\put(20,60){\circle*{3}}
\put(20,40){\circle*{3}}
\put(30,70){\circle*{3}}
\put(30,30){\circle*{3}}
\put(40,70){\circle*{3}}
\put(40,30){\circle*{3}}
\put(50,70){\circle*{3}}
\put(50,30){\circle*{3}}
\put(10,10){$\al$}
\put(30,10){$\phi (\xi_1)$}

\put(90,70){\line(1,0){65}}
\put(120,70){\line(-2,-1){35}}
\put(115,30){\line(1,0){40}}
\put(95,30){\line(1,0){10}}
\put(120,30){\line(-2,-1){35}}
\put(120,70){\vector(0,-1){55}}

{\thicklines
\put(110,50){\line(1,-1){20}}
\put(110,50){\line(1,1){10}}
\put(120,60){\line(2,-3){20}}
\put(130,70){\line(1,-2){20}}
\put(140,70){\line(1,-2){10}}

\put(110,65){\line(-1,-3){5}}
\put(105,50){\line(1,-5){5}}
\put(100,30){\line(0,1){5}}
\put(100,40){\line(0,1){5}}
\put(100,50){\line(0,1){5}}
\put(100,65){\line(0,1){5}}
\put(90,15){\line(0,1){40}}
}
\put(100,70){\circle*{3}}
\put(100,30){\circle*{3}}
\put(110,66){\circle*{3}}
\put(110,25){\circle*{3}}
\put(90,55){\circle*{3}}
\put(90,15){\circle*{3}}
\put(120,60){\circle*{3}}
\put(120,40){\circle*{3}}
\put(130,70){\circle*{3}}
\put(130,30){\circle*{3}}
\put(140,70){\circle*{3}}
\put(140,30){\circle*{3}}
\put(150,70){\circle*{3}}
\put(150,30){\circle*{3}}
\put(110,10){$\al$}
\put(130,10){$\phi (\eta_1)$}

\put(195,70){\line(1,0){80}}
\put(220,15){\line(1,0){55}}
\put(205,15){\line(1,0){10}}
\put(230,15){\line(-2,-1){35}}
\put(230,70){\line(-2,-1){35}}
\put(230,70){\vector(0,-1){65}}
\put(220,0){$\al$}
\put(240,-5){$\phi (\sm_1)$}

\put(230,20){\circle*{3}}
\put(230,35){\circle*{3}}
\put(230,50){\circle*{3}}
\put(230,65){\circle*{3}}
\put(240,70){\circle*{3}}
\put(250,70){\circle*{3}}
\put(260,70){\circle*{3}}
\put(240,15){\circle*{3}}
\put(250,15){\circle*{3}}
\put(260,15){\circle*{3}}
\put(200,0){\circle*{3}}
\put(200,55){\circle*{3}}
\put(205,15){\circle*{3}}
\put(205,70){\circle*{3}}

{\thicklines
\put(240,70){\line(-2,-1){10}}
\put(225,63){\line(-2,-1){10}}
\put(210,45){\line(0,1){10}}
\put(220,40){\line(2,-1){10}}
\put(250,70){\line(-1,-1){20}}
\put(210,40){\line(2,1){20}}
\put(210,40){\line(0,-1){10}}
\put(240,15){\line(-2,1){30}}
\put(250,15){\line(-1,1){20}}
\put(260,15){\line(0,1){55}}
\put(200,0){\line(0,1){55}}
\put(205,15){\line(0,1){5}}
\put(205,25){\line(0,1){5}}
\put(205,35){\line(0,1){5}}
\put(205,45){\line(0,1){5}}
\put(205,65){\line(0,1){5}}
}

\put(320,70){\line(1,0){80}}
\put(340,20){\line(1,0){60}}
\put(324,20){\line(1,0){8}}
\put(350,70){\line(-2,-1){35}}
\put(350,20){\line(-2,-1){35}}
\put(350,70){\vector(0,-1){65}}
\put(340,0){$\al$}
\put(365,0){$\phi (\tau_1)$}

\put(350,25){\circle*{3}}
\put(350,45){\circle*{3}}
\put(350,65){\circle*{3}}
\put(360,70){\circle*{3}}
\put(375,70){\circle*{3}}
\put(385,70){\circle*{3}}
\put(360,20){\circle*{3}}
\put(375,20){\circle*{3}}
\put(385,20){\circle*{3}}
\put(330,70){\circle*{3}}
\put(330,20){\circle*{3}}
\put(335,64){\circle*{3}}
\put(335,14){\circle*{3}}
\put(320,55){\circle*{3}}
\put(320,5){\circle*{3}}

{\thicklines
\put(350,45){\line(-1,1){10}}
\put(350,45){\line(-1,-1){10}}
\put(350,45){\line(1,1){25}}
\put(350,45){\line(1,-1){25}}
\put(340,55){\line(1,1){10}}
\put(340,35){\line(1,-1){10}}
\put(360,70){\line(-2,-1){10}}
\put(360,20){\line(-2,1){10}}

\put(385,20){\line(0,1){50}}
\put(335,13){\line(0,1){50}}
\put(320,6){\line(0,1){50}}
\put(330,20){\line(0,1){5}}
\put(330,30){\line(0,1){5}}
\put(330,40){\line(0,1){5}}
\put(330,50){\line(0,1){5}}
\put(330,65){\line(0,1){5}}
}

\put(150,-20){{\bf Fig.~5:} Three-page singular tangles}
\end{picture}
\vspace{1.0cm}

As for singular tangles from Subsection~3.1, isotopy classes
 of three-page tangles in the layer $\{0\leq z\leq 1\}$ form a
semigroup.
Each three-page tangle can be encoded by a word in the alphabet
 $\A=\lbrace a_i,b_i,c_i,d_i,x_i \vl i \in \Z_3\rbrace$ (Fig.~3)
in the same way as in Subsection~2.3.
A three-page tangle is called \emph{almost balanced}, if the
corresponding word in the alphabet $\A$ is 1-balanced and
2-balanced (see Subsection~2.4).
Note that for any $i$-balanced 3-page tangle all strings in the
band $P_i$ can be assumed vertical.
By $BT$ we denote the semigroup of almost balanced 3-page tangles.
Define the map $\ph:ST\to BT$ on the generators as follows, $k\in \N$
(Fig.~5):
 $$(24)\; \phi(\xi_k)=d_2^k c_2 b_2^{k-1},
 \phi(\eta_k)=d_2^{k-1} a_2 b_2^k,
 \phi(\sm_k)=d_2^{k-1} b_1 d_2 d_1 b_2^k,
 \phi(\si)=d_2^k b_1 b_2 d_1 b_2^{k-1},
 \phi(\tau_k)=d_2^k x_2 b_2^k.$$
Each tangle goes to the corresponding three-page embedding plus
vertical intervals.
The following Lemma is proved in the same way as \cite[Lemma~3]{Dy3}.

\begin{lem}
The map $\ph:ST\to BT$ is well-defined isomorphism of semigroups.
\end{lem}
\begin{proof}
First let us check that isotopy equivalent singular tangles go to
isotopy equivalent 3-page singular tangles under the map $\ph$.
Actually, by definition, singular tangles from semigroups $ST$ and
 $BT$ are considered up to isotopy in the layer $\{0<z<1\}$.
So the injectivity of the map $\ph$ follows.
Now we construct the inverse map $\psi: BT\to ST$.
Let us associate with each almost balanced 3-page tangle $\Ga\in BT$
the singular tangle $\psi(\Ga)\in ST$ given by the following diagram.
According to the almost balance we consider that all segments of
$\Ga$ lying in the pages $P_1,P_2$ are vertical.
Deleting all these vertical segments from $\Ga$ we obtain some
graph-tangle $\psi(\Ga)$ in the sense of Subsection~3.1.
Clearly, the composition $\psi\circ\ph: ST\to ST$ is
 identical on the generators (Fig.~5).
So, the maps $\ph,\psi$ are mutually inverse.
\end{proof}

\subsection{Classification of three-page singular tangles}
By $\ph(11)-\ph(23)$ we denote the relations between words
in the alphabet $\A$ which are obtained from the relations (11)-(23)
 of the semigroup $ST$ under the isomorphism $\ph$ (Lemma~2).
The following Lemma will be proved in Subsection~4.3.

\begin{lem}
Relations (1)-(10) follow from relations $\ph(11)-\ph(23)$ of the
semigroup $SK$.
\end{lem}

Theorem~2 is a special case of the following classification theorem
 for three-page singular tangles, which we prove by analogy with
Theorem~1 of \cite{Dy3}.

\begin{theo}
The semigroup of the isotopy classes of 3-page singular tangles is
isomorphic to the semigroup $SK$.

\end{theo}
\begin{proof}
As it was already mentioned in Subsection~3.2 with each three-page
singular tangle it is possible to associate a word in the alphabet
$\A$ and hence an element of the semigroup $SK$. Conversely, each
element of the semigroup $SK$ can be completed to some three-page
tangle, if we add three families of parallel segments on each page
$P_i$, $i\in \Z_3$. For example, the three-page tangles
corresponding to the following elements from $SK$ are depicted on
Fig.~3: $d_2 c_2$, $a_2 b_2$, $b_1 d_2 d_1 b_2$, $d_2 x_2 b_2$.
Relations (1)-(10) of the semigroup $SK$ can be easily performed
 by an ambient isotopy in the layer $\{0<z<1\}$.
Hence it remains to prove that each isotopy of three-page singular
 tangle can be decomposed into ''elementary isotopies" corresponding to
relations (1)-(10) of $SK$.
It suffices to do this for almost balanced three-page tangles.
Really, for a 3-page tangle given by a word $w$
 let $n_1$ and $m_1$ (respectively, $n_2$ and $m_2$)
 be the maximal numbers of points on the semilines
 $P_1\cap \{z=0\}$ and $P_1\cap \{z=1\}$ (respectively, on the
 semilines $P_2\cap \{z=0\}$ and $P_2\cap \{z=1\}$) that are
 connected by arcs with points on the segment $I$.
For letter $a_0$ one gets $n_1=n_2=0,m_1=m_2=1$ (Fig.~2a).
For an almost balanced 3-page graph-tangle we have
 $n_1=m_1=n_2=m_2=0$.
Then the word $b_1^{n_2} d_2^{n_1} w b_2^{m_1} d_1^{m_2}$
 is almost balanced.
Because of invertibility of generators $b_i,d_i$ such a transformation
sends equivalent words to equivalent.
By Lemma~2 we can associate a singular tangle in the sense of
Subsection~3.1 with each almost balanced word.
For such tangles each isotopy is already decomposed into the
elementary isotopies corresponding to relations $\ph(11)-\ph(23)$
(Lemmas~1 and 2). So Lemma~3 finishes the proof of Theorem~4.
\end{proof}

\subsection{Proof of Theorem~3}
We will identify an arbitrary three-page singular tangle with
the corresponding element of the semigroup $SK$. A 3-page tangle
is called \emph{knot-like} if it contains a singular knot near
the axis $\al$ and the rest of it consists only of vertical
segments. Evidently, the knot-like tangles correspond to
balanced words from $W_b\subset W$.

\begin{lem}
An element $w\in SK$ defines a knot-like tangle if and only if
$w$ is a central element in the semigroup $SK$.
\end{lem}
\begin{proof}
The part ''only if" is geometrically evident: a singular knot can
be moved by an isotopy to any place of a given tangle, i.e., a
knot-like element commutes with any other by Theorem~2.
Let $w$ be a central element in $SK$. Then for each $k\in \N$ we have
 $b_i^k d_i^k w = w b_i^k d_i^k$.
Denote by $m$ (respectively by $n$) the number of arcs of the
3-page tangle $w$ that go out in the page $P_{i-1}$ to the left
(respectively, to the right) boundary.
Then for sufficiently large $k$ the number of arcs of the 3-page
tangle $b_i^k d_i^k w$ that go out in the page $P_{i-1}$ to the
left boundary is equal to $k$, and for the tangle $w b_i^k d_i^k$
is equal to $m+k-n$, i.e. $m=n$. Hence for sufficiently large $l$
and any $j=i-1\in \Z_3$ the word $a_0^l a_1^l w c_1^l c_0^l$ is
$j$-balanced, so it is balanced.
Since $w$ is central, the word $w a_0^l a_1^l c_1^l c_0^l$ is also
balanced. Then it is geometrically obvious that the element $w$ defines
 a knot-like tangle.
\end{proof}
Theorem~3 follows from Lemma~4.
\subsection{Classification of knuckle 4-valent graphs}
The isotopy classification problem of such graphs was considered
 in the paper \cite{Ka}.
Dynnikov's method gives a possibility of solving it analogously
 to the case of singular knots.
Let us introduce a semigroup $FG$ with the same generators
 and relations as the semigroup $SK$, only we change
 relation (6) to the following:
$(6') \quad x_i (d_{i+1} d_i d_{i-1})=x_i.$
Then the semigroup $FG$ has 15 generators and 84 defining relations.

\begin{theo}
The center of the semigroup $FG$ classifies all non-oriented
 knuckle 4-valent graphs in $\R^3$ up to an ambient isotopy.
\end{theo}

Theorem~5 is proved analogously to Theorems~1--3 with the change of
relation (6) to $(6')$ and relation (23) in Lemma~1 to the following
 $(23')\; \sm_k \tau_k=\tau_k$.

\section{Proof of Propositions~1 and 2}
In Claim~1 we obtain new word equivalences from relations (1)-(10)
of the semigroup $SK$.
In Subsection~4.2, using Claims~1--3 we prove Lemma~5 about
decomposition of any $i$-balanced word.
Lemma~5 and Claim~6 reduce the infinite number of relations
 $\ph(11)-\ph(23)$ to the finite number of relations (1)-(10).
The proof of Lemma~3 finishes in Subsection~4.3 using Claims~5
and 6. All relations will be obtained in a formal way, but they
 have a geometric interpretation (Fig.~3).
\subsection{Corollaries of relations (1)-(10)}
\begin{claim}
Equivalences (1)-(10) imply the following (we suppose that
$i\in \Z_3$ and
$w_i\in \B_{i}=$ $\{  a_i, \; b_i,$ $ \; c_i, \; d_i, \; b_{i-1} b_i
d_{i-1}, \; b_{i-1} d_i d_{i-1} \; \}$):
$$\begin{array}{ll}
(25) \; b_i \s d_{i+1} d_{i-1}, \mbox{ or } b_{i+1} \s d_{i-1} d_i,
     \; b_{i-1} \s d_i d_{i+1}, & \mbox{ or }
(25) \; b_0 \s d_1 d_2, \; b_1\s d_2 d_0, \; b_2 \s d_0 d_1; \\

(26) \; d_i \s b_{i-1} b_{i+1}, \mbox{ or } d_{i-1} \s b_{i+1} b_i,
     \; d_{i+1} \s b_i b_{i-1},  & \mbox{ or }
(26) \; d_0 \s b_2 b_1, \; d_1\s b_0 b_2, \; d_2 \s b_1 b_0; \\

(27) \; d_{i+1} b_{i-1} \s b_{i-1} d_{i+1} t_i,
     \; b_{i+1} d_{i-1} \s t_i d_{i-1} b_{i+1},  & \mbox{ where }
        t_i=b_{i+1} d_{i-1} d_{i+1} b_{i-1};
\end{array} $$
$$\begin{array}{ll}
(28) \; a_i \s a_{i-1} b_{i+1}, \; c_i \s d_{i+1} c_{i-1}; &
(29) \; a_i b_i \s a_{i-1} d_{i-1},  \; d_i c_i \s b_{i-1} c_{i-1}; \\

(30) \; b_i \s a_i b_i c_i, \; d_i \s a_i d_i c_i; & \\

(31) \; b_{i-1} x_{i+1} d_{i-1} \s x_i; &
(32) \; b_{i} x_i d_{i} \s d_{i+1} x_{i+1} b_{i+1};
\end{array} $$
$$\begin{array}{ll}
(33) \; ( d_{i} c_{i} ) w_{i+1} \s w_{i+1} ( d_{i} c_{i} ); &
(34) \; ( b_{i} c_{i} ) w_{i-1} \s w_{i-1} ( b_{i} c_{i} ); \\

(35) \; ( a_{i} b_{i} ) w_{i+1} \s w_{i+1} ( a_{i} b_{i} ); &
(36) \; ( a_{i} d_{i} ) w_{i-1} \s w_{i-1} ( a_{i} d_{i} );  \\

(37) \; t_{i} w_{i} \s w_{i} t_{i}, \; t_{i}' w_{i} \s w_{i} t_{i}',
 \mbox{ where }  &
 t_i=b_{i+1} d_{i-1} d_{i+1} b_{i-1}, \;
 t_i'=d_{i-1} b_{i+1} b_{i-1} d_{i+1};
\end{array} $$
$$\begin{array}{ll}
(38) \; (d_i x_i b_i) w_{i+1} \s w_{i+1} (d_i x_i b_i); &
(39) \; (b_i x_i d_i) w_{i-1} \s w_{i-1} (b_i x_i d_i); \\

(40) \; d_{i+1} b_{i-1} w_i d_{i-1} b_{i+1} \s
        b_{i-1} d_{i+1} w_i b_{i+1} d_{i-1}; &
\end{array} $$
$$\begin{array}{ll}
(41) \; b_{i-1}^2 a_i d_{i-1}^2 \s
  ( b_{i-1} a_i d_{i-1} ) d_i^2 ( b_{i-1} b_i d_{i-1} ) b_i; &

(42) \;  b_{i-1}^2 c_i d_{i-1}^2 \s
  d_i ( b_{i-1} d_i d_{i-1} ) b_i^2 ( b_{i-1} c_i d_{i-1} ); \\

(43) \; b_{i-1}^2 b_i d_{i-1}^2 \s
  ( b_{i-1} b_i d_{i-1} ) d_i^2 ( b_{i-1} b_i d_{i-1} ) b_i; &

(44) \; b_{i-1}^2 d_i d_{i-1}^2 \s
  d_i ( b_{i-1} d_i d_{i-1} ) b_i^2 ( b_{i-1} d_i d_{i-1} );
\end{array}$$

$$(45) \; b_{i-1}^2 x_i d_{i-1}^2 \s
  ( b_{i-1}^2 b_i d_{i-1}^2 ) ( b_{i-1} d_i d_{i-1} )
  d_i^2 x_i b_i^2 ( b_{i-1} b_i d_{i-1} ) ( b_{i-1}^2 d_i d_{i-1}^2 ).$$
\end{claim}
\begin{proof}
Note that equivalences (25)-(27) follow easily from (3)-(4).
By (4) we have $d_i\s b_i^{-1}$,
$b_{i-1} b_i d_{i-1} \s (b_{i-1} d_i d_{i-1})^{-1}$, and
$t_i'\s t_i^{-1}$.
Then (35), (37) and (38) follow from (8), (9) and (10) respectively.
The other equivalences will be verified step by step, using those already
checked. Recall that $i\in \Z_3=\{0,1,2\}$, i.e.
 in $\Z_3$ we have $(i+1)+1=i-1$ and $(i-1)-1=i-1$.

$$\begin{array}{ll}
(28) \quad
  a_{i-1} b_{i+1}
   \st{(1)}{\s}
  (a_i d_{i+1}) b_{i+1}
   \st{(4)}{\s} a_i, &

  d_{i+1} c_{i-1}
   \st{(1)}{\s}
  d_{i+1} (b_{i+1} c_i)
   \st{(4)}{\s}
  c_i; \\

(29) \quad
  a_i b_i
   \st{(28)}{\s}
  (a_{i-1} b_{i+1}) b_i
   \st{(26)}{\s}
  a_{i-1} d_{i-1},  &

  d_i c_i
   \st{(26)}{\s}
  (b_{i-1} b_{i+1}) c_i
   \st{(1)}{\s}
  b_{i-1} c_{i-1}; \\

(30) \quad
  a_i b_i c_i
   \st{(29)}{\s}
  a_{i} (d_{i+1} c_{i+1})
   \st{(1)}{\s}
  a_{i-1} c_{i+1}
   \st{(1)}{\s}
  b_{i},  &

  a_i d_i c_i
   \st{(29)}{\s}
  a_{i} (b_{i-1} c_{i-1})
   \st{(28)}{\s}
  a_{i+1} c_{i-1}
   \st{(1)}{\s}
  d_{i}; \\

(31) \;
  b_{i-1} x_{i+1} d_{i-1}
   \st{(2)}{\s}
  b_{i-1} (d_{i-1} x_{i} b_{i-1}) d_{i-1}
   \st{(4)}{\s}
  x_i; \\ 

(32) \;
  b_{i} x_{i} d_{i}
   \st{(31)}{\s}
  b_i (b_{i-1} x_{i+1} d_{i-1}) d_{i}
   \st{(25), \ (26)}{\s}
  d_{i+1} x_{i+1} b_{i+1}.
\end{array}$$

Below in the proof of (33) we firstly commute $b_{i+1}$ with
 $d_i c_i$, and then we use this relation to commute $a_{i+1}$
 with $d_i c_i$.
$$\begin{array}{l}
(33b) \quad
  b_{i+1} (d_i c_{i})
   \st{(30)}{\s}
  (a_{i+1} b_{i+1} c_{i+1}) (d_i c_i)
   \st{(7)}{\s}
  a_{i+1} b_{i+1} (d_i c_i) c_{i+1}
   \st{(26)}{\s}
  a_{i+1} b_{i+1} (b_{i-1} b_{i+1}) c_i c_{i+1}
   \st{(1)}{\s} \\ \st{(1)}{\s}
  (a_{i+1} b_{i+1}) b_{i-1} c_{i-1} c_{i+1}
   \st{(8)}{\s}
  b_{i-1} c_{i-1} (a_{i+1} b_{i+1}) c_{i+1}
   \st{(30)}{\s}
  b_{i-1} c_{i-1} b_{i+1}
   \st{(1)}{\s}
  b_{i-1} (b_{i+1} c_i) b_{i+1}
   \st{(26)}{\s}
  (d_{i} c_i) b_{i+1}; \\ \\

(33a) \quad
  a_{i+1} (d_i c_{i})
   \st{(1)}{\s}
  (a_{i-1} d_{i}) (d_i c_{i})
   \st{(26)}{\s}
  a_{i-1} (b_{i-1} b_{i+1}) (d_i c_{i})
   \st{(33b)}{\s}
  a_{i-1} b_{i-1} (d_i c_i) b_{i+1}
   \st{(35)}{\s} \\ \qquad \st{(35)}{\s}
  (d_{i} c_{i}) (a_{i-1} b_{i-1}) b_{i+1}
   \st{(26)}{\s}
  (d_{i} c_{i}) (a_{i-1} d_{i})
   \st{(1)}{\s}
  (d_{i} c_{i}) a_{i+1}.
\end{array} $$

The other equivalences in (33) follow from (33a), (33b) and (7).
Equivalences (34), (36),(39) follow respectively from (29) and (33),
(29) and (35), (32) and (38).
The last calculations are trivial:
$$\begin{array}{l}
(40) \;
  d_{i+1} b_{i-1} w_i d_{i-1} b_{i+1}
   \st{(27)}{\s}
  (b_{i-1} d_{i+1} t_i) w_i d_{i-1} b_{i+1}
   \st{(39)}{\s}
  b_{i-1} d_{i+1} (w_i t_i) d_{i-1} b_{i+1}
   \st{(27)}{\s}
  b_{i-1} d_{i+1} w_i b_{i+1} d_{i-1}; \\ \\

(41) \;
 b_{i-1}^2 a_i d_{i-1}^2
  \st{(4)}{\s}
 b_{i-1}^2 a_i (d_i b_i) d_{i-1}^2
  \st{(36)}{\s}
 b_{i-1} (a_i d_{i}) (b_{i-1} b_{i}) d_{i-1}^2
  \st{(26)}{\s}
 b_{i-1} a_{i} d_{i} (b_{i-1} b_{i}) d_{i-1} (b_{i+1} b_i)
  \st{(37)}{\s} \\ \qquad \st{(37)}{\s}
 b_{i-1} a_i b_{i+1} (d_{i} b_{i-1} b_i d_{i-1}) b_i
  \st{(25)}{\s}
 (b_{i-1} a_i d_{i-1}) d_{i}^2 (b_{i-1} b_i d_{i-1}) b_i; \\ \\

(42) \;
 b_{i-1}^2 b_i d_{i-1}^2
  \st{(4)}{\s}
 b_{i-1} (b_i d_i) b_{i-1} b_i d_{i-1}^2
  \st{(26)}{\s}
 b_{i-1} b_i (d_i b_{i-1} b_i d_{i-1}) (b_{i+1} b_i)
  \st{(37)}{\s} \\ \qquad \st{(37)}{\s}
 b_{i-1} b_i b_{i+1} (d_i b_{i-1} b_i d_{i-1}) b_i
  \st{(25)}{\s}
 (b_{i-1} b_i d_{i-1}) d_i^2 (b_{i-1} b_i d_{i-1}) b_i;
\end{array}$$
$$\begin{array}{l}
(43) \;
 b_{i-1}^2 c_i d_{i-1}^2
  \st{(4)}{\s}
 b_{i-1}^2 (d_i b_i) c_i d_{i-1}^2
  \st{(34)}{\s}
 b_{i-1}^2 d_{i} d_{i-1} (b_{i} c_i) d_{i-1}
  \st{(25)}{\s}
 (d_{i} d_{i+1}) (b_{i-1} d_i d_{i-1} b_{i}) c_i d_{i-1}
  \st{(37)}{\s} \\ \qquad \st{(37)}{\s}
 d_{i} (b_{i-1} d_i d_{i-1} b_{i}) d_{i+1} c_i d_{i-1}
  \st{(26)}{\s}
 d_{i} (b_{i-1} d_i d_{i-1}) b_{i}^2 (b_{i-1} c_i d_{i-1}); \\ \\

(44) \;
 b_{i-1}^2 d_i d_{i-1}^2
  \st{(4)}{\s}
 b_{i-1}^2 d_i d_{i-1} (b_i d_i) d_{i-1}
  \st{(25)}{\s}
 (d_i d_{i+1}) b_{i-1} d_i d_{i-1} b_i d_{i} d_{i-1}
  \st{(37)}{\s} \\ \qquad \st{(37)}{\s}
 d_i (b_{i-1} d_i d_{i-1} b_i) d_{i+1} d_{i} d_{i-1}
  \st{(26)}{\s}
 d_i (b_{i-1} d_i d_{i-1}) b_i^2 (b_{i-1} d_{i} d_{i-1});
\end{array}$$
$$\begin{array}{l}
(45) \;
 b_{i-1}^2 x_i d_{i-1}^2
  \st{(4)}{\s}
 b_{i-1}^2 (b_i d_i) x_i (b_i b_{i+1}^2 d_{i+1}^2 d_i) d_{i-1}^2
  \st{(10)}{\s}
 b_{i-1}^2 b_i b_{i+1}^2 (d_i x_i b_i) d_{i+1}^2 d_i d_{i-1}^2
  \st{(25),(26)}{\s} \\ \st{(25),(26)}{\s}
 b_{i-1}^2 b_i (d_{i-1} d_i)^2 d_i x_i b_i (b_i b_{i-1})^2 d_i d_{i-1}^2
  \st{(4)}{\s}
 ( b_{i-1}^2 b_i d_{i-1}^2 ) ( b_{i-1} d_i d_{i-1} )
   d_i^2 x_i b_i^2 ( b_{i-1} b_i d_{i-1} ) ( b_{i-1}^2 d_i d_{i-1}^2 ).
 \end{array}$$
\end{proof}
\subsection{Decomposition of $i$-balanced words}
\begin{claim}
For each $i\in \Z_3$ any $i$-balanced word is equivalent to some
 $i$-balanced word that contains only the following letters:
 $a_i$, $b_i$, $c_i$, $d_i$, $x_i$, $b_{i-1}$, $d_{i-1}$.
\end{claim}

\begin{proof}
Using the following substitutions, we can eliminate the other
letters:
$$\left \lbrace \begin{array}{l}
  a_{i-1} \st{(1)}{\s} a_i d_{i+1}, \quad
  a_{i+1} \st{(28)}{\s} a_i b_{i-1}, \quad
  c_{i-1} \st{(1)}{\s} b_{i+1} c_i, \quad
  c_{i+1} \st{(28)}{\s} d_{i-1} c_i, \\
  b_{i+1} \st{(25)}{\s} d_{i-1} d_i, \quad
  d_{i+1} \st{(26)}{\s} b_i b_{i-1}, \quad
  x_{i-1} \st{(2)}{\s} d_i x_{i+1} b_i, \quad
  x_{i+1} \st{(2)}{\s} d_{i-1} x_{i} b_{i-1}. \quad
  \end{array} \right.$$
\end{proof}
Fix an index $i\in\Z_3$. Let $w$ be an $i$-balanced word on the
letters $a_i$, $b_i$, $c_i$, $d_i$, $x_{i}$, $b_{i-1}$, $d_{i-1}$.
Consider the substitution
 $\mu: a_i,b_i,c_i,d_i,x_{m,i}\to \bu$,
 $b_{i-1}\to ($, $d_{i-1}\to )$.
Denote by $\mu(w)$ the resulting \emph{encoding} consisting of
brackets and bullets. Because the given word $w$ is
$i$-balanced, so the encoding $\mu(w)$ (without bullets)
is a balanced bracket expression. For each place
$k$ denote by $dif(k)$ the difference
between the number of left and right brackets in a subword of
$\mu(w)$, ending at this place. The maximum of
$dif(k)$ for all $k$ we call the \emph{depth}
of $K$: $d(w)$. For example, the word $w=b_{i-1}^2 a_i d_{i-1}^2$
has the encoding $\mu(w)=((\bu))$ and the depth $d(w)=2$.

By the \emph{star of depth} $k$ we call an encoding of the type
$(^k \bu )^k$ which has $k$ couples of brackets. The bullet is a
star of depth 0. If for a word $w$ its encoding $\mu(w)$
decomposes into several stars, then $w$ is called
\emph{star decomposable}.
In this case the depth $d(w)$ is maximal among the depths of all
stars participating in the decomposition.
\begin{claim}
Every $i$-balanced word $w$ is equivalent to some star
decomposable word $w'$ of the same depth $d(w')=d(w)$.
\end{claim}
\begin{proof}
Consider the beginning of the encoding $\mu(w)$. After several
initial left brackets $\mu(w)$ contains either a right bracket
or a bullet. In the first case we delete a couple of brackets
$()$ by the rule:
 $b_{i-1} d_{i-1}\st{(2)}{\s} \emptyset$.
Hence we can suppose that the next simbol after the sequence of
$k$ left brackets is a bullet. Because $\mu(w)$ is balanced,
then after this bullet it may be the sequence of $j$,
$0 \leq j \leq k$, right brackets. If $j<k$, then insert
into $w$ after the last right bracket the following subword
$d_{i-1}^{k-j} b_{i-1}^{k-j}\st{(2)}{\s} \es$.
This operation does not change the depth $d(w)$. Then in the
resulting word $w_1$ the encoding $\mu(w_1)$ contains a star of
depth $k$ at the beginning. Continuing this process, after a finite
number of steps we get a star decomposable word $w_N$ of the same
depth $d(w_N)=d(w)$.
\end{proof}
For any letter $s$ denote by $s^\prime$ the word
$b_{i-1} s d_{i-1}$, for example, $a_i'=b_{i-1} a_i d_{i-1}$.
\begin{claim}
Each star decomposable word $w$ is equivalent to a word
decomposed into the following $i$-balanced subwords:
 $ \{ a_i, \; b_i, \; c_i, \; d_i, \; x_{i}, \;
 a_i', \; b_i', \; c_i', \; d_i', \; x_{i}' \}.$
\end{claim}
\begin{proof} We use induction on the depth $d(w)$.
The case $d(w)=1$ is trivial. Let $\mu(w)$ contain a star
of depth $k\geq 2$. Apply one of the following
transformations to every such star.
$$\left\{ \begin{array}{l}
 u=b_{i-1}^2 a_i d_{i-1}^2 \st{(46a)}{\s} a_i' d_i^2 b_i' b_i = v,
 \mbox{ i.e. } \mu(u)=((\bu)) \to \mu(v)=(\bu)\bu\bu(\bu)\bu; \\

 u=b_{i-1}^2 b_i d_{i-1}^2 \st{(46b)}{\s} b_i' d_i^2 b_i' b_i = v,
 \mbox{ i.e. } \mu(u)=((\bu)) \to \mu(v)=(\bu)\bu\bu(\bu)\bu; \\

 u=b_{i-1}^2 c_i d_{i-1}^2 \st{(46c)}{\s} d_i d_i' b_i^2 c_i' = v,
 \mbox{ i.e. } \mu(u)=((\bu)) \to \mu(v)=\bu(\bu)\bu\bu(\bu); \\

 u=b_{i-1}^2 d_i d_{i-1}^2 \st{(46c)}{\s} d_i d_i' b_i^2 d_i' = v,
 \mbox{ i.e. } \mu(u)=((\bu)) \to \mu(v)=\bu(\bu)\bu\bu(\bu);
\end{array} \right.$$

$$ \left\{ \begin{array}{l}
    u=b_{i-1}^2 x_{i} d_{i-1}^2
   \st{(46x)}{\s}
   ( b_{i-1}^2 b_i d_{i-1}^2 ) ( b_{i-1} d_i d_{i-1} )
   d_i^2 x_{i} b_i^2 ( b_{i-1} b_i d_{i-1} )
( b_{i-1}^2 d_i d_{i-1}^2 )
   \st{(46b,46d)}{\s} \\
  \st{(46b,46d)}{\s}
   (b_i' d_i^2 b_i' b_i) d_i' d_i^2 x_{i} b_i^2 b_i'
(d_i d_i' b_i^2 d_i') = v,
   \mbox{ i.e.  } \\
  \mu(u)=((\bu)) \to
   \mu(v)=(\bu) \bu \bu (\bu) \bu (\bu) \bu \bu \bu \bu \bu
(\bu) \bu (\bu) \bu \bu (\bu).
\end{array} \right.  $$
We get  a word $w_1\s w$, equivalent to a star decomposed
word $w_2$ of depth $d(w_2)=d(w)-1$ according to Claim~3.
This finishes the induction step.
\end{proof}
\begin{lem}
For all $i\in \Z_3$ each $i$-balanced word is equivalent
 to some word which can be decomposed into the $i$-balanced words
 of the set
 $\B_i=\{ a_i, b_i, c_i, d_i, x_i, b_{i-1} b_i d_{i-1},
  b_{i-1} d_i d_{i-1} \}.$
\end{lem}
\begin{proof}
By Claims~3 and 4 it remains to eliminate only the following words:

$$\begin{array}{l}
   a'_i = b_{i-1} a_i d_{i-1}
    \st{(25)}{\s}
   (d_i d_{i+1}) a_i d_{i-1}
    \st{(4)}{\s}
   d_i d_{i+1} a_i (b_i d_i) d_{i-1}
    \st{(35)}{\s}
   d_i (a_i b_i) d_{i+1} d_i d_{i-1}
    \st{(26)}{\s}
   d_i a_i b_i^2 (b_{i-1} d_i d_{i-1}); \\

   c'_i = b_{i-1} c_i d_{i-1}
    \st{(26)}{\s}
   b_{i-1} c_i (b_{i+1} b_i)
    \st{(4)}{\s}
   b_{i-1} (b_i d_i c_i) b_{i+1} b_i
    \st{(33)}{\s}
   b_{i-1} b_i b_{i+1} (d_i c_i) b_i
    \st{(35)}{\s}
   (b_{i-1} b_i d_{i-1}) d_i^2 c_i b_i;
\end{array}$$
$$\left\lbrace  \begin{array}{l}
   x'_i = b_{i-1} x_i d_{i-1}
    \st{(4)}{\s}
   b_{i-1} (b_i b_{i+1} d_{i+1} d_i) x_i (b_i d_i) d_{i-1}
    \st{(38)}{\s}
   b_{i-1} b_i b_{i+1} (d_i x_i b_i) d_{i+1} d_i d_{i-1}
    \st{(25)}{\s} \\ \st{(25)}{\s}
   (b_{i-1} b_i d_{i-1}) d_i^2 x_i b_i d_{i+1} d_i d_{i-1}
    \st{(26)}{\s}
   (b_{i-1} b_i d_{i-1}) d_i^2 x_i b_i^2 (b_{i-1} d_i d_{i-1}).
  \end{array} \right. $$
\end{proof}
Denote by $(33')-(40')$ relations (33)-(40), assuming
that $w_i\in W_i$.
\begin{claim}
Generalized equivalences $(33')-(40')$ hold for
 arbitrary $i$-balanced words $w_i\in W_i$.
\end{claim}
\begin{proof}
By Lemma~5 every $i$-balanced word of $W_i$ can be decomposed
into the elementary words from $\B_i$. Since equivalences (33)-(40)
hold for words from $\B_i$, they also hold for words from $W_i$.
Note that we get an infinite number of new equivalences
 $(33')-(40')$.
\end{proof}

\subsection{Deduction of relations $\ph(11)-\ph(23)$ from relations
(1)-(10) of the semigroup $SK$}
For each $l\geq 1$ denote by $u_l$ the symbol from the set of
generators $\{ \xi_l, \eta_l, \sm_l, \sm_l^{-1}, \tau_l \}$
of the semigroup $ST$. Define the \emph{shift maps}
$\te_k: ST\to ST$ and
 $\rho_k: BT\to BT$ by $\te_k(u_l)=u_{k+l}$ and
 $\rho_k(w)=d_2^k w b_2^k$.
Evidently, the shift map $\te_k: ST\to ST$ is a well-defined
 homomorphism.
Really, each relation from (11)-(23) for $k>1$ is obtained
 from the corresponding relation for $k=1$ by the shift map
 $\te_{k-1}$.
For example, relation $\xi_k \xi_l=\xi_{l+2} \xi_k$ is
 obtained from $\xi_1 \xi_{l-k+1}=\xi_{l-k+3} \xi_1$ by
the shift map $\te_{k-1}$.
By (2) the shift map $\rho_k$ sends equivalent words to
 equivalent ones, i.e. $\rho_k$ is a homomorphism.
Moreover, the following diagram
$$\begin{CD}
ST      @>{\te_k}>>  ST \\
@V{\ph}VV              @VV{\ph}V \\
BT      @>{\rho_k}>> ST
\end{CD}$$
is commutative, that implies Claim~6.
\begin{claim}
For each $k\in\N$ relations $\ph(11)-\ph(23)$ can be obtained from
relations $\ph(11)-\ph(23)$ for $k=1$ on using the equivalences
$b_2 d_2\s 1\s d_2 b_2$ of (2).
\qed
\end{claim}
\begin{proof}[{\bf Proof of Lemma~3}]
Here we deduce relations $\ph(11)-\ph(23)$ of the semigroup $BT$
among the words in the alphabet $\A$ from relations (1)-(10) and
(25)-(40). We use generalized equivalences $(33')-(40')$ from
Claim~5. By the star $\star$ we denote the following images of the
map $\ph$ for $k=1$:
 $$(24')\; \ph(\xi_1)=d_2 c_2, \;
           \ph(\eta_1)=a_2 b_2, \;
           \ph(\sm_1)=b_1 d_2 d_1 b_2, \;
           \ph(\sm_1^{-1})=d_2 b_1 b_2 d_1, \;
           \ph(\tau_1)=d_2 x_2 b_2.$$
Then by (24) we have $\ph(u_l)=d_2^{l-1} \star b_2^{l-1}$.
Note that the words $d_2^{l} \star b_2^{l}\in W$ are
 1-balanced and 2-balanced (Fig.~5).
Then relations $\ph(11)-\ph(14)$ are proved by the same scheme:


$$\begin{array}{l}
(11) \;
  \ph(\xi_1 u_l)
   \st{(24)}{=}
  (d_2 c_2) (d_2^{l-1} \star b_2^{l-1})
   \st{(4)}{\s}
  d_2^2 (b_2 c_2) (d_2^{l-1} \star b_2^{l-1})
   \st{(34')}{\s}
  d_2^2 (d_2^{l-1} \star b_2^{l-1}) (b_2 c_2)

   \st{(4)}{\s}
  \ph(u_{l+2} \xi_1); \\

(12) \;
  \ph(\eta_1 u_l)
   \st{(24)}{=}
  (a_2 b_2) (d_2^{l-1} \star b_2^{l-1})
   \st{(4)}{\s}
  (a_2 d_2) (d_2^{l-3} \star b_2^{l-3}) b_2^2
   \st{(36')}{\s}
  (d_2^{l-3} \star b_2^{l-3}) (a_2 d_2) b_2^2
   \st{(4)}{\s}
  \ph(u_{l-2} \eta_1); \\

(13) \;
  \ph(\sm_1 u_l)
   \st{(24)}{=}
  (b_1 d_2 d_1 b_2) (d_2^{l-1} \star b_2^{l-1})
   \st{(4)}{\s}
  (b_1 d_2 d_1) (d_2^{l-2} \star b_2^{l-1})
   \st{(25),(2)}{\s}
  d_2^2 (b_2 d_0 d_2 b_0) (d_2^{l-3} \star b_2^{l-3}) b_2^2
   \st{(37')}{\s} \\ \qquad \st{(37')}{\s}
  d_2^2 (d_2^{l-3} \star b_2^{l-3}) (b_2 d_0 d_2 b_0) b_2^2
   \st{(4),(26)}{\s}
  (d_2^{l-1} \star b_2^{l-2}) (b_2 b_1) d_2 d_1 b_2
   \st{(24)}{=}
  \ph(u_{l} \sm_1); \\

(14) \;
  \ph(\tau_1 u_l)
   \st{(24)}{=}
  (d_2 x_2 b_2) (d_2^{l-1} \star b_2^{l-1})
   \st{(4)}{\s}
  d_2^2 (b_2 x_2 d_2) (d_2^{l-3} \star b_2^{l-3}) b_2^2
   \st{(39')}{\s}
  d_2^2 (d_2^{l-3} \star b_2^{l-3}) (b_2 x_2 d_2) b_2
  \st{(4)}{\s} \\ \qquad \st{(4)}{\s}
  (d_2^{l-1} \star b_2^{l-1}) (d_2 x_2 b_2)
   \st{(24)}{=}
  \ph(u_{l} \tau_1).
\end{array} $$

The remaining calculations are straightforward:


$$(15) \;
  \ph(\eta_{2} \xi_1)
   \st{(24)}{=}
  (d_2 a_2 b_2^{2}) (d_2 c_2)
   \st{(4)}{\s}
  d_2 (a_2 b_2 c_2)
   \st{(30)}{\s}
  d_2 b_2
   \st{(4)}{\s}
  1
   \st{(4)}{\s}
  d_2 b_2
   \st{(30)}{\s}
  (a_2 d_2 c_2) b_2
   \st{(4)}{\s}
  \ph(\eta_1 \xi_{2}); $$

$$(16) \; \left \lbrace \begin{array}{l}
  \ph(\eta_{3} \sm_{3} \xi_1)
   \st{(24),(4)}{\s}
  d_2^{2} a_2 b_2 (b_2 b_1) d_2 d_1 (b_2 c_2)
   \st{(26)}{\s}
  d_2^{2} (a_2 b_2) d_0 d_2 d_1 (b_2 c_2)
   \st{(35)}{\s}
  d_2^{2} d_0 (a_2 b_2) d_2 d_1 (b_2 c_2)
   \st{(34)}{\s} \\ \st{(34)}{\s}
  d_2 (d_2 d_0) (a_2 b_2) d_2 (b_2 c_2) d_1
   \st{(4)}{\s}
  d_2 (d_2 d_0) (a_2 b_2 c_2) d_1
   \st{(25)}{\s}
  d_2 b_1 (a_2 b_2 c_2) d_1
   \st{(30)}{\s}
  d_2 b_1 b_2 d_1 b_2
   \st{(24)}{=}
  \ph(\sm_1^{-1}), \\ \\

  \ph(\eta_1 \sm_{2} \xi_{3})
   \st{(4)}{\s}
  a_2 b_1 d_2 (d_1 d_2) c_2 b_2^{2}
   \st{(25)}{\s}
  a_2 b_1 d_2 (b_0 c_2) b_2^{2}
   \st{(1)}{\s}
  (a_0 d_1) b_1 d_2 c_1 b_2^{2}
   \st{(4)}{\s}
  a_0 d_2 (b_1 d_1) c_1 b_2^{2}
   \st{(33)}{\s} \\ \st{(33)}{\s}
  a_0 d_2 b_1 b_2 (d_1 c_1) b_2
   \st{(37)}{\s}
  (d_2 b_1 b_2 d_1) (a_0 c_1) b_2
   \st{(1)}{\s}
  d_2 b_1 b_2 d_1 d_2 b_2
   \st{(4)}{\s}
  d_2 b_1 b_2 d_1
   \st{(24)}{=}
  \ph(\sm_1^{-1});
 \end{array} \right. $$

$$(17) \;
  \ph(\eta_{3} \tau_{3} \xi_1)
   \st{(24),(4)}{\s}
  d_2^{2} a_2 b_2 x_2 b_2 c_2
   \st{(5)}{\s}
  d_2^{2} (b_2 x_2 b_2)
   \st{(4)}{\s}
  \ph(\tau_1)
   \st{(4)}{\s}
  (d_2 x_2 d_2) b_2^{2}
   \st{(5)}{\s}
  a_2 d_2 x_2 d_2 c_2 b_2^{2}
   \st{(4)}{\s}
  \ph(\eta_{1} \tau_{2} \xi_{3});$$

$$(18) \left \lbrace \begin{array}{l}
  \ph(\sm_1 \xi_1)
   \st{(24)}{=}
  (b_1 d_2 d_1 b_2) (a_2 b_2)
   \st{(4)}{\s}
  b_1 d_2 (d_1 c_2)
   \st{(28)}{\s}
  b_1 (d_2 c_0)
   \st{(28)}{\s}
  b_1 c_1
   \st{(25)}{\s}
  (d_2 d_0) c_1
   \st{(28)}{\s}
  a_2 b_2
   \st{(24)}{=}
  \ph(\xi_1); \\

  \ph(\eta_1 \sm_1)
   \st{(4)}{\s}
  a_2 (b_2 b_1) d_2 d_1 b_2
   \st{(26)}{\s}
  (a_2 d_0) d_2 d_1 b_2
   \st{(1)}{\s}
  (a_1 d_2) d_1 b_2
   \st{(1)}{\s}
  (a_0 d_1) b_2
   \st{(1)}{\s}
  \ph(\eta_k);
 \end{array} \right. $$


$$(19) \;
  \ph(\sm_1 \sm_1^{-1})
   \st{(24)}{=}
  (b_1 d_2 d_1 b_2) (d_2 b_1 b_2 d_1)
   \st{(4)}{\s}
  (b_1 d_2 d_1) (b_1 b_2 d_1)
   \st{(4)}{\s}
  1
   \st{(4)}{\s}
  (d_2 b_1 b_2 d_1) (b_1 d_2 d_1 b_2)
   \st{(24)}{=}
  \ph(\sm_1^{-1} \sm_1);$$


$$(20) \; \left\lbrace \begin{array}{l}
  \ph(\sm_{2} \sm_1 \sm_{2})
   \st{(4)}{\s}
  d_2 b_1 d_2 d_1 b_2^2 b_1 d_2^2 d_1 b_2^{2}
   \st{(26)}{\s}
  d_2 (b_1 d_2 d_1 b_2) d_0 d_2^2 d_1 b_2^{2}
   \st{(37)}{\s}
  d_2 d_0 (b_1 d_2 d_1 b_2) d_2^2 d_1 b_2^{2}
   \st{(25),(4)}{\s} \\ \st{(25),(4)}{\s}
  b_1^2 d_2 (d_1 d_2) d_1 b_2^{2}
   \st{(26)}{\s}
  b_1^2 d_2 (d_1 d_2) d_1 (d_0 d_1) b_2
   \st{(25)}{\s}
  b_1^2 d_2 b_0 d_1 (d_0 d_1) b_2
   \st{(4)}{\s}
  b_1^2 (d_2 b_0 d_1 d_0 b_2) d_2 d_1 b_2
   \st{(40)}{\s} \\ \st{(40)}{\s}
  b_1^2 (b_0 d_2 d_1 b_2 d_0) d_2 d_1 b_2
   \st{(25)}{\s}
  b_1^2 (d_1 d_2) d_2 d_1 b_2 (b_2 b_1) d_2 d_1 b_2
   \st{(26)}{\s}
  b_1 d_2^2 d_1 b_2 (b_2 b_1) d_2 d_1 b_2
   \st{(4)}{\s}
  \ph(\sm_1 \sm_{2} \sm_1);
 \end{array} \right.$$

$$(21) \; \left\lbrace \begin{array}{l}
  \ph(\tau_{2} \sm_1 \sm_{2})
   \st{(4)}{\s}
  d_2^{2} x_2 b_2^2 b_1 d_2^2 d_1 b_2^{2}
   \st{(26)}{\s}
  d_2 (d_2 x_2 b_2) d_0 d_2^2 d_1 b_2^{2}
   \st{(38)}{\s}
  d_2 d_0 (d_2 x_2 b_2) d_2^2 d_1 b_2^{2}
   \st{(4)}{\s} \\ \st{(4)}{\s}
  (d_2 d_0) (d_2 x_2 d_2) d_1 b_2^{2}
   \st{(25)}{\s}
  b_1 (d_2 x_2 d_2) d_1 b_2^{2}
   \st{(4)}{\s}
  b_1 d_2^2 (b_2 x_2 d_2) d_1 b_2^{2}
   \st{(39)}{\s}
  b_1^2 d_2^2 d_1 (b_2 x_2 d_2) b_2^{2}
   \st{(4)}{\s}
  \ph(\sm_1 \sm_{2} \tau_1).
 \end{array} \right.$$

$$(22) \left\lbrace \begin{array}{l}
  \ph(\tau_{1} \sm_{2} \sm_{1})
   \st{(4)}{\s}
  d_2 x_2 (b_1 d_2 d_1 b_2) (b_2 b_1) d_2 d_1 b_2
   \st{(26)}{\s}
  d_2 x_2 (b_1 d_2 d_1 b_2) d_0 d_2 d_1 b_2
   \st{(37),(4)}{\s} \\ \st{(37),(4)}{\s}
  d_2 (d_0 b_0) x_2 d_0 (b_1 d_2 d_1^2) b_2
   \st{(31)}{\s}
  d_2 d_0 x_1 (b_1 d_2 d_1^2) b_2
   \st{(4)}{\s}
  d_2 d_0 b_1 (d_1 x_1 b_1) d_2 d_1^2 b_2
   \st{(38)}{\s} \\ \st{(38)}{\s}
  d_2 d_0 b_1 d_2 (d_1 x_1 b_1) d_1^2  b_2
   \st{(2),(4)}{\s}
  d_2 d_0 b_1 d_2 d_1 (b_0 x_2 d_0) d_1 b_2
   \st{(25),(4)}{\s}
  d_2 d_0 (b_1 d_2 d_1 b_2) d_2 b_0 x_2 b_2^{2}
   \st{(37)}{\s} \\ \st{(37)}{\s}
  d_2 (b_1 d_2 d_1 b_2) d_0 d_2 b_0 x_2 b_2^{2}
   \st{(25),(26)}{\s}
  d_2 b_1 d_2 d_1 b_2 (b_2 b_1) d_2 (d_1 d_2) x_2 b_2^{2}
   \st{(4)}{\s}
  \ph(\sm_{2} \sm_{1} \tau_{2});
 \end{array} \right.$$

$$(23) \;
  \ph(\sm_{1} \tau_1)
   \st{(4)}{\s}
  b_1 d_2 d_1 x_2 b_2
   \st{(25)}{\s}
  (d_2 d_0) d_2 d_1 x_2 b_2
   \st{(6)}{\s}
  d_2 x_2 (d_0 d_2 d_1) b_2
   \st{(26)}{\s}
  d_2 x_2 (b_2 b_1) d_2 d_1 b_2
   \st{(4)}{\s}
  \ph(\tau_1 \sm_{1}).$$

\end{proof}

\end{document}